# QUANTILE ESTIMATION WITH ADAPTIVE IMPORTANCE SAMPLING

By Daniel Egloff and Markus Leippold

*QuantCatalyst and University of Zurich*

We introduce new quantile estimators with adaptive importance sampling. The adaptive estimators are based on weighted samples that are neither independent nor identically distributed. Using a new law of iterated logarithm for martingales, we prove the convergence of the adaptive quantile estimators for general distributions with nonunique quantiles thereby extending the work of Feldman and Tucker [*Ann. Math. Statist.* **37** (1996) 451–457]. We illustrate the algorithm with an example from credit portfolio risk analysis.

**1. Introduction.** We introduce a new sample-based quantile estimators with adaptive importance sampling. Importance sampling is a widely used technique for variance reduction to improve the statistical efficiency of Monte Carlo simulations. It reduces the number of samples required for a given level of accuracy. The basic idea is to change the sampling distribution so that a greater concentration of samples are generated in a region of the sample space which has a dominant impact on the calculations. The change of distribution is then compensated by weighting the samples using the Radon–Nikodym derivative of the original measure with respect to the new measure. However, in a multivariate setting, it is far from obvious how such a change of measure should be obtained.

Given its importance for practical applications, especially for risk management in the finance industry, the literature on sample-based quantile estimation with variance reduction is rather sparse.[1] The focus of variance reduction schemes is almost exclusively geared towards the estimation of expected values. The reason might lie in the additional intricateness that sample-based quantile estimators exhibit. The quantile function, viewed as a map on the space of distribution functions, generally fails to be differentiable in the sense of Hadamard. For certain distributions, the quantile may be nonunique. If the lower and upper quantiles of a random variable $Y$ for a probability level $\alpha \in (0,1)$, defined as

$$q_\alpha(Y) = \inf\{y \mid \mathbb{P}(Y \le y) \ge \alpha\},$$
$$q^\alpha(Y) = \sup\{y \mid \mathbb{P}(Y \le y) \le \alpha\},$$

are distinct, then the ordinary quantile estimator $Y_{\lfloor n\alpha \rfloor, n}$ based on the order statistics $Y_{1,n}, \ldots, Y_{n,n}$ of samples $Y_1, \ldots, Y_n$ is not consistent anymore. For the special case of independent and identically distributed (i.i.d.) samples $Y_1, \ldots, Y_n$, Feldman and Tucker [12] prove that $Y_{\lfloor n\alpha \rfloor, n}$ oscillates across the interval $[q_\alpha(Y), q^\alpha(Y)]$ infinitely often. Using the classical law of iterated logarithm for sequences of i.i.d. random variables, they also show that consistency can be retained for the modified estimator $Y_{\nu(n),n}$ if the function $\nu(n) \in \mathbb{N}$ satisfies

$$(1+k)\sqrt{2n \log\log n} \le \lfloor n\alpha \rfloor - \nu(n) \le K n^{1/2+\gamma} \tag{1.1}$$

for some positive constants $\gamma, k, K$.

[1]An early contribution is the control variate approach of [18]. Later work includes [4, 16] and [19].

For the estimation of expected values with importance sampling, a common procedure is to apply the change of measure suggested by a large deviation upper bound. Although this approach often leads to an asymptotically optimal sampling estimator, it can also fail completely, as shown in Glasserman and Wang [14].

Our method for quantile estimation does not rely on large deviation principles. Instead, it is adaptive. Adaptive algorithms, but only for expected values and not for quantiles, are introduced in the work of [2] and [3]. They apply the truncated Robbins–Monro algorithm of Chen, Guo and Gao in [7] for pricing financial options under different assumptions on the underlying process. Robbins–Monro methods and stochastic approximation date back to the historical work of Robbins and Monro [30] and Kiefer and Wolfowitz [21]. See also [23] and the more recent references [5, 24] and [25].

Using an adaptive strategy to obtain a quantile estimator means that every new sample is used to improve the parameters of the importance sampling density. Therefore, we cannot rely on the results of Feldman and Tucker [12]. Our quantile estimators, derived from weighted samples, are neither independent nor identically distributed. However, we derive a new law of iterated logarithms for martingales which allows us to prove the convergence of the adaptive quantile estimators for distributions with nondifferentiable and nonunique quantiles without requiring the i.i.d. assumption thereby extending the result of Feldman and Tucker.

Our paper is structured as follows. In Section 2, we present the general setup and we introduce the notation. Section 3 gives a brief review of adaptive importance sampling for estimating the mean. In Section 4, we start with the discussion of the metric structure underlying our adaptive algorithm. We then derive two theorems, Theorems 4.1 and 4.2, which extend Feldman and Tucker [12]. The proof of the theorems build on a new result for the law of iterated logarithms for martingales which we present in Theorem 4.4. Finally, in Section 5 we provide an application of our new quantile estimator which we borrow from credit risk management. All proofs are delegated to the Appendix.

**2. Setup and notation.** Let $\varphi_\theta(x)$ be a probability density depending on a parameter $\theta$ of a random variable $X$ relative to some reference measure $\lambda$, defined on a measurable space $(\boldsymbol{\mathcal{X}}, \mathcal{F})$ with a countably generated $\sigma$-field $\mathcal{F}$. We assume that the parameters $\theta$ take their values in a metric space $(\boldsymbol{\Theta}, d)$ for some fixed metric $d$ and equip it with the Borel $\sigma$-algebra $\mathcal{B}(\boldsymbol{\Theta})$. For now, we do not have to further distinguish the parameter space $\boldsymbol{\Theta}$ and the set of densities $\{\varphi_\theta(x) \mid \theta \in \boldsymbol{\Theta}\}$. For the expectation, under the measure $\varphi_\theta \, d\lambda$, we write

$$\mathbb{E}_\theta[f(X)] = \int_{\boldsymbol{\mathcal{X}}} f(x)\varphi_\theta(x)\, d\lambda(x) \tag{2.1}$$

and we define for all $p$, $1 \le p \le \infty$,

$$\mathcal{L}_p(\theta) = \{f : \boldsymbol{\mathcal{X}} \to \mathbb{R} \mid f \text{ is } \mathcal{F}\text{-measurable}, \|f\|_{\theta,p}^p = \mathbb{E}_\theta[|f(X)|^p] < \infty\}.$$

Let $\varphi_{\theta_0}(x)$ be our reference or sample density. We assume that all densities in $\boldsymbol{\Theta}$ are absolutely continuous with respect to the reference density $\varphi_{\theta_0}(x)$ and we denote by

$$w_x(\theta) = \frac{\varphi_{\theta_0}(x)}{\varphi_\theta(x)} \tag{2.2}$$

the likelihood ratio or Radon–Nikodym derivative. In particular, $w_x : x \mapsto w_x(\theta)$ is measurable for all $\theta \in \boldsymbol{\Theta}$. If $f \in \mathcal{L}_1(\theta_0)$, we have

$$\mathbb{E}_\theta[w_X(\theta) f(X)] = \mathbb{E}_{\theta_0}[f(X)] \qquad \forall \theta \in \boldsymbol{\Theta}. \tag{2.3}$$

For $p \ge 1$, we introduce the (weighted) moments

$$m_{f,p}(\theta) = \mathbb{E}_\theta[|w_X(\theta) f(X)|^p] = \mathbb{E}_{\theta_0}[w_X(\theta)^{p-1}|f(X)|^p]. \tag{2.4}$$

We use the abbreviation $m_f(\theta) = m_{f,2}(\theta)$ for the second moment and

$$\sigma_f^2(\theta) = \mathrm{Var}_\theta[w_X(\theta) f(X)] = m_f(\theta) - \mathbb{E}_{\theta_0}[f(X)]^2 \tag{2.5}$$

for the variance.

**3. Review: Adaptive importance sampling for estimation of means.** Before we derive our adaptive quantile estimators, we start this section with a brief review of adaptive importance sampling for the estimation of means. Consider a function $f \in \mathcal{L}_1(\theta_0)$. Static importance sampling estimates the expectation $\mathbb{E}_{\theta_0}[f(X)]$ by the weighted sample average

$$\hat{e}_{\mathrm{s}}(n,f) = \frac{1}{n}\sum_{i=1}^{n} w_{X_i}(\theta) f(X_i) \tag{3.1}$$

with $X_i \sim \varphi_\theta \, d\lambda$ i.i.d. The usual error estimates based on the central limit theorem indicate that the most advantageous choice for $\theta$ would be the variance minimizer

$$\theta^* = \arg\min_{\theta \in \boldsymbol{\Theta}} \sigma_f^2(\theta) = \arg\min_{\theta \in \boldsymbol{\Theta}} m_f(\theta). \tag{3.2}$$

Unfortunately, in most cases (3.2) cannot be solved explicitly. An alternative to the approach on the basis of large deviation upper bounds is an adaptive strategy. The solution $\theta^*$ is estimated by a sequence $(\theta_n)_{n\geq 0}$, generated, for instance, by a stochastic approximation algorithm of Kiefer–Wolfowitz or Robbins–Monro type. Replacing the fixed parameter $\theta$ in (3.1) by the sequentially generated parameters $(\theta_n)_{n\geq 0}$ leads to the adaptive importance sampling estimator

$$\hat{e}_{\mathrm{a}}(n,f) = \frac{1}{n}\sum_{i=1}^{n} w_{X_i}(\theta_{i-1}) f(X_i), \tag{3.3}$$

where $X_n \sim \varphi_{\theta_{n-1}}(x)\, d\lambda(x)$ is simulated from the importance sampling distribution determined from the parameter $\theta_{n-1}$. In contrast to static importance sampling, the random variables $w_{X_n}(\theta_{n-1})$ and $f(X_n)$ in (3.3) are neither independent nor identically distributed. However, we still obtain a martingale.

LEMMA 3.1. *Let $\theta_n$ be a sequence of parameters and $X_n \sim \varphi_{\theta_{n-1}}\, d\lambda$. Define $\mathcal{F}_n = \sigma(\theta_0,\dots,\theta_n, X_1,\dots,X_n)$. Then, for $f \in \mathcal{L}_1(\theta_0)$,*

$$M_n = \sum_{i=1}^{n} (w_{X_i}(\theta_{i-1}) f(X_i) - \mathbb{E}_{\theta_0}[f(X)]) \tag{3.4}$$

*is a martingale with respect to the filtration $\mathbb{F} = (\mathcal{F}_n)_{n\geq 0}$.*

PROOF. If $(Z_n)_{n\geq 0}$ is a sequence of integrable random variables, then

$$M_n = \sum_{i=1}^{n} (Z_i - \mathbb{E}[Z_i \mid Z_1,\dots,Z_{i-1}]) \tag{3.5}$$

is a martingale. The integrability of $w_{X_i}(\theta_{i-1})f(X_i)$ and the martingale property for (3.4) follow from

$$\mathbb{E}[w_{X_i}(\theta_{i-1})f(X_i) \mid \mathcal{F}_{i-1}] = \mathbb{E}_{\theta_{i-1}}[w_{X_i}(\theta_{i-1})f(X_i)] = \mathbb{E}_{\theta_0}[f(X)], \tag{3.6}$$

where the second equality is a consequence of (2.3). □

A strong law of large numbers and a central limit theorem for (3.3) has been obtained in [2] by applying classical martingale convergence results for which we refer to [17] and [15]. For a proof of the theorem below, we refer to [2].

Theorem 3.1. *Let $\theta_n$, $X_n$, and $\mathbb{F} = (\mathcal{F}_n)_{n\geq 0}$ be as in Lemma 3.1. Assume that $\theta_n \to \theta^* \in \mathbf{\Theta}$ converges almost surely and that there exists $a > 1$ such that for all $\theta \in \mathbf{\Theta}$*

$$\mathbb{E}_\theta[|w_X(\theta)f(X)|^{2a}] < \infty, \tag{3.7}$$

*the function $m_{f,2a} : \theta \mapsto m_{f,2a}(\theta)$ is continuous in $\theta^*$, and*

$$\mathbb{E}[m_{f,2a}(\theta_n)] < \infty \qquad \forall n \geq 0. \tag{3.8}$$

*Then*

$$\lim_{n\to\infty} \frac{1}{n}\sum_{i=1}^{n} w_{X_i}(\theta_{i-1})f(X_i) = \mathbb{E}_{\theta_0}[f(X)] \qquad \text{almost surely}, \tag{3.9}$$

*and*

$$\sqrt{n}\left(\frac{1}{n}\sum_{i=1}^{n} w_{X_i}(\theta_{i-1})f(X_i) - \mathbb{E}_{\theta_0}[f(X)]\right) \xrightarrow{d} N(0, \sigma_f^2(\theta^*)), \tag{3.10}$$

*where $\xrightarrow{d}$ denotes convergence in distribution.*

**4. Adaptive importance sampling for quantile estimation.** Having reviewed the estimation of the mean with adaptive importance sampling in the previous section, we introduce now the metric structure and the algorithm that underlies our new adaptive quantile estimation.

4.1. *Riemannian structure for parameter tuning.* The procedure to estimate the variance optimal parameter (3.2) crucially depends on the metric structure of the parameter space $\mathbf{\Theta}$. The metric is not only important if it comes to the actual numerical implementation, but is also material to determine existence and uniqueness of a solution.

Let the parameter space $\mathbf{\Theta}$ be a smooth manifold. It is known that the canonical metric on a family of densities $\{\varphi_\theta(x) \mid \theta \in \mathbf{\Theta}\}$ is induced by the Riemannian structure given by the Fisher information metric

$$g = \mathbb{E}_\theta[dl_X \otimes dl_X]. \tag{4.1}$$

Here, $l_x(\theta) = \log \varphi_\theta(x)$ is the log-likelihood function with differential

$$dl_x : \mathbf{\Theta} \to T^*\mathbf{\Theta}, \tag{4.2}$$

considered as a one-form on $\mathbf{\Theta}$ and with $T^*\mathbf{\Theta}$, the co-tangent space of the manifold $\mathbf{\Theta}$. In particular, (4.1) defines a nondegenerate symmetric bilinear form on the tangent space $T\mathbf{\Theta}$, hence a Riemannian metric.[2]

Having equipped the parameter space $\mathbf{\Theta}$ with a Riemannian metric, we can formulate the first order condition for (3.2) in terms of the Riemannian gradient $\nabla$ as

$$\nabla m_f(\theta) = 0. \tag{4.3}$$

Under suitable assumptions on $X$ and the likelihood ratio $w_x(\theta)$, we can exchange integration and differentiation to arrive at

$$\nabla m_f(\theta) = -\mathbb{E}_{\theta_0}[f(X)^2 w_X(\theta) \nabla l_X(\theta)] = -\mathbb{E}_\theta[f(X)^2 w_X(\theta)^2 \nabla l_X(\theta)] \tag{4.4}$$

with $\nabla l_x(\theta)$ the Riemannian gradient of the log-likelihood. To approximate a solution $\nabla m_f(\theta^*) = 0$, we can now use the representation (4.4) and a stochastic approximation scheme

$$\theta_{n+1} = \theta_n + \gamma_{n+1} H(X_{n+1}, \theta_n), \qquad X_{n+1} \sim \varphi_{\theta_n}\, d\lambda, \tag{4.5}$$

with average descent direction

$$H(X, \theta) = -f(X)^2 w_X(\theta)^2 \nabla l_X(\theta). \tag{4.6}$$

In this paper, we want to keep the focus on adaptive importance sampling for quantile estimation and we therefore restrict ourselves to vector spaces. An example with a flat metric is provided by the Gaussian densities

$$\mathbf{\Theta} = \{N(\theta, \Sigma) \mid \theta \in \mathbb{R}^k\} \tag{4.7}$$

with fixed covariance structure $\Sigma$.[3] The first and second order differentials of the likelihood $l_x(\theta)$ are

$$dl_x(\theta) = \Sigma^{-1}(x - \theta), \qquad d^2 l_x(\theta) = -\Sigma^{-1}. \tag{4.10}$$

---

[2] For the basic concepts of Riemannian geometry, we refer to [20, 22] and the references therein. The usage of the Riemannian metric based on the Fisher information goes back to [28].

[3] A well-known example of a nonflat Riemannian structure on a space of distributions is the Fisher information metric of a location scale family of densities

$$\mathbf{\Theta} = \left\{ \varphi_{(\mu,\sigma)}(x) = \frac{1}{\sigma}\varphi\left(\frac{x-\mu}{\sigma}\right) \middle| (\mu, \sigma) \in \mathbb{R} \times \mathbb{R}^+ \right\}. \tag{4.8}$$

A second example is given by the space of all multivariate normal distributions

$$\mathbf{\Theta} = \{N(\theta, \Sigma) \mid \theta \in \mathbb{R}^k, \Sigma \in \mathcal{S}_+(k)\} \tag{4.9}$$

which is not flat anymore.

Hence, the Fisher metric on $\mathbf{\Theta}$ is

$$g_\Sigma(u,v) = -\mathbb{E}_\theta[d^2 l_X(\theta)(u,v)] = u^\top \Sigma^{-1} v. \tag{4.11}$$

Because

$$g_\Sigma(\nabla l_x(\theta), u) = \nabla l_x(\theta)^\top \Sigma^{-1} u = dl_x(\theta)(u) = (x-\theta)^\top \Sigma^{-1} u, \tag{4.12}$$

the gradient of the likelihood with respect to the metric (4.11) is

$$\nabla l_x(\theta) = (x-\theta). \tag{4.13}$$

Note that the gradient $\nabla$ of the Fisher metric defers by a factor of $\Sigma^{-1}$ from the gradient induced by the standard Euclidian metric.

4.2. *Parameter tuning with adaptive truncation.* In practical applications, the parameter space $\mathbf{\Theta}$ is often noncompact or it is difficult to a priori identify a bounded region to which the optimal parameter must belong. We therefore suppose that the parameter sequence $(\theta_n)_{n\geq 0}$ is generated by an algorithm that enforces recurrence and boundedness of the sequence $\theta_n$ by adaptive truncation. A series of work [6–10] shows that stochastic approximation algorithms with adaptive truncation behave numerically more smoothly and converge under weaker hypotheses. No restrictive conditions on the mean field or a-priori boundedness assumptions have to be imposed. We follow Andrieu, Moulines and Priouret [1] who analyze the convergence of stochastic approximation algorithms with more flexible truncation schemes and Markov state-dependent noise. To specify the algorithm, we let $(\mathcal{K}_j)_{j\in\mathbb{N}}$ be an increasing compact covering of $\mathbf{\Theta}$ satisfying

$$\mathbf{\Theta} = \bigcup_{j=1}^{\infty} \mathcal{K}_j, \mathcal{K}_j \subset \operatorname{int}(\mathcal{K}_{j+1}) \tag{4.14}$$

and

$$\boldsymbol{\gamma} = (\gamma_n)_{n\in\mathbb{N}}, \qquad \boldsymbol{\epsilon} = (\epsilon_n)_{n\in\mathbb{N}}, \tag{4.15}$$

two monotonically decreasing sequences. We introduce the counting variables

$$(\kappa_n, \nu_n, \zeta_n)_{n\in\mathbb{N}} \in \mathbb{N}\times\mathbb{N}\times\mathbb{N}, \tag{4.16}$$

where $\kappa_n$ records the active truncation set in the compact covering, $\nu_n$ counts the number of iterations since the last re-initialization (truncation) and $\zeta_n$ is the index in the sequences $\boldsymbol{\gamma}$, $\boldsymbol{\epsilon}$ introduced in (4.15). If $\nu_n \neq 0$, the algorithm operates in the active truncation set $\mathcal{K}_{\kappa_n}$ so that

$$\theta_j \in \mathcal{K}_{\kappa_n} \qquad \forall j \leq n \text{ with } \nu_j \neq 0. \tag{4.17}$$

If $\nu_n = 0$, the update at iteration $n$ caused a jump outside of the active truncation set $\mathcal{K}_{\kappa_n}$ and triggers a re-initialization at the next iteration $n+1$. We assume that a stochastic vector field is generated from a measurable map

$$H : \boldsymbol{\mathcal{X}} \times \boldsymbol{\Theta} \to \boldsymbol{\Theta}, \tag{4.18}$$

where $\boldsymbol{\mathcal{X}}$ and $\boldsymbol{\Theta}$ are both equipped with countably generated $\sigma$-fields $\mathcal{B}(\boldsymbol{\mathcal{X}})$ and $\mathcal{B}(\boldsymbol{\Theta})$, respectively. We also suppose that $\boldsymbol{\Theta}$ is an open subset of some Euclidian vector space.

To handle jumps outside the parameter space $\boldsymbol{\Theta}$, we introduce an isolated point $\theta_c$ taking the role of a cemetery point. Let $\bar{\boldsymbol{\Theta}} = \boldsymbol{\Theta} \cup \{\theta_c\}$. For an arbitrary $\gamma \geq 0$, we define a kernel $Q_\gamma$ on $\boldsymbol{\mathcal{X}} \times \bar{\boldsymbol{\Theta}}$ by

$$\begin{aligned} Q_\gamma(x,\theta; A \times B) = & \int_A P_\theta(dy) \mathbf{1}_{\{\theta + \gamma H(y,\theta) \in B\}} \\ & + \mathbf{1}_{\{\theta_c \in B\}} \int_A P_\theta(dy) \mathbf{1}_{\{\theta+\gamma H(y,\theta) \notin B\}}, \end{aligned} \tag{4.19}$$

where $(x,\theta) \in \boldsymbol{\mathcal{X}} \times \boldsymbol{\Theta}$ and $A \in \mathcal{B}(\boldsymbol{\mathcal{X}})$, $B \in \mathcal{B}(\bar{\boldsymbol{\Theta}})$, and $P_\theta(dx)$ is a measure on $\mathcal{X}$ parameterized by $\theta$. For a sequence of step sizes $\boldsymbol{\gamma}$, we define the process $(X_n, \theta_n)$ by

$$(X_{n+1}, \theta_{n+1}) \sim Q_{\gamma_{n+1}}(X_n, \theta_n; \cdot) \tag{4.20}$$

unless $\theta_n = \theta_c$, in which case we stop the process and set $\theta_{n+1} = \theta_c$, $X_{n+1} = X_n$. The law of the nonhomogeneous Markov process (4.20) with initial conditions $(x,\theta)$, represented on the product space $(\boldsymbol{\mathcal{X}} \times \bar{\boldsymbol{\Theta}})^{\mathbb{N}}$, is denoted by $\mathbb{P}^{\boldsymbol{\gamma}}_{x,\theta}$. Let $\mathcal{X}_0 \subset \boldsymbol{\mathcal{X}}$ be a compact subset

$$\Phi : \boldsymbol{\mathcal{X}} \times \boldsymbol{\Theta} \to \mathcal{X}_0 \times \mathcal{K}_0, \tag{4.21}$$

be a measurable map and $\phi : \mathbb{N} \to \mathbb{Z}$ with $\phi(n) > -n$.

ALGORITHM 4.1. The stochastic approximation algorithm with adaptive truncation is the homogeneous Markov chain defined by the following transition law from step $n$ to $n+1$:

(i) If $\nu_n = 0$, then we perform a reset operation which starts in $\mathcal{X}_0 \times \mathcal{K}_0$ and draws

$$(X_{n+1}, \theta_{n+1}) \sim Q_{\gamma_{\zeta_n}}(\Phi(X_n, \theta_n); dx \times d\theta).$$

Otherwise, we simulate

$$(X_{n+1}, \theta_{n+1}) \sim Q_{\gamma_{\zeta_n}}(X_n, \theta_n; dx \times d\theta).$$

(ii) If $\|\theta_{n+1} - \theta_n\| \le \epsilon_{\zeta_n}$ and $\theta_{n+1} \in \mathcal{K}_{\kappa_n}$, then we update

$$\nu_{n+1} = \nu_n + 1, \qquad \zeta_{n+1} = \zeta_n + 1, \qquad \kappa_{n+1} = \kappa_n;$$

else we prepare for a reset operation in the next iteration by setting

$$\nu_{n+1} = 0, \qquad \zeta_{n+1} = \phi(\zeta_n), \qquad \kappa_{n+1} = \kappa_n + 1.$$

The convergence of Algorithm 4.1 under suitable conditions on the measure $P_\theta(dx)$, the mean field $h$ defined as

$$h(\theta) = \int_{\boldsymbol{\mathcal{X}}} H(x,\theta) P_\theta(dx) \tag{4.22}$$

and the sequences $\boldsymbol{\gamma}$, $\boldsymbol{\epsilon}$ are established in [1].[4]

4.3. *Quantile estimation.* After having introduced the metric structure and the parameter tuning in the previous sections, we can now turn our focus to the estimation of quantiles of a real-valued random variable

$$Y = \Psi(X), \qquad \Psi : \mathcal{X} \to \mathbb{R},$$

defined in terms of an $\mathcal{F}$-measurable function $\Psi$. We denote by

$$q_\alpha = q_\alpha(Y) = \inf\{y \mid \mathbb{P}(Y \le y) \ge \alpha\}, \qquad 0 < \alpha < 1,$$

the lower $\alpha$-quantile of $Y$. Furthermore, let $(\theta_n)_{n \ge 0}$ be a sequence of tuning parameters. In favor of a more compact notation, we introduce the abbreviations

$$w_n = w_{X_n}(\theta_{n-1}), \qquad Y_n = \Psi(X_n), \qquad n \ge 1. \tag{4.23}$$

We recall that, under the assumptions of Theorem 3.1, the weights $w_n$ satisfy

$$\mathbb{E}[w_n \mid \mathcal{F}_{n-1}] = \mathbb{E}_{\theta_{n-1}}[w_n] = 1, \qquad \frac{1}{n}\sum_{i=1}^n w_i \to 1 \qquad \text{almost surely.} \tag{4.24}$$

We first consider generalizations of the empirical distribution function to weighted samples. Because the sum of the weights $\sum_{i=1}^n w_i$ is not necessarily normalized to one, we introduce the renormalized weighted empirical distribution function

$$F_{n,w}(y) = \frac{1}{\sum_{i=1}^n w_i} \sum_{i=1}^n w_i \mathbf{1}_{\{Y_i \le y\}} \tag{4.25}$$

[4] In fact, [1] treat the more general case of state dependent noise where the measure $P_\theta(dy)$ takes the form of a Markov kernel $P_\theta(x, dy)$.

and set

$$F_{n,w,\nu}(y) = \frac{1}{\nu(n)} F_{n,w}(y), \tag{4.26}$$

where $\nu : \mathbb{N} \to \mathbb{R}^+$ is a normalization function, which we determine later to prevent the oscillation of the quantile estimators. We can use the increasing function $F_{n,w,\nu}$ to define the quantile estimator

$$q_{n,w,\nu}(\alpha) = F^{\leftarrow}_{n,w,\nu}(\alpha) = \inf\{y \mid F_{n,w,\nu}(y) \geq \alpha\}, \tag{4.27}$$

where $F^{\leftarrow}_{n,w,\nu}$ is the generalized inverse of $F_{n,w,\nu}$. Besides the re-normalized weighted empirical distribution function (4.25), there are alternative ways to generalize the empirical distribution function to weighted samples. For example,

$$F^{l}_{n,w,\nu}(y) = \frac{1}{\nu(n)} \sum_{i=1}^{n} w_i \mathbf{1}_{\{Y_i \leq y\}} \tag{4.28}$$

puts the emphasis on the left tail of the distribution. However, if the concern is on the right tail, then

$$\begin{aligned} F^{r}_{n,w,\nu}(y) &= \frac{1}{\nu(n)} \sum_{i=1}^{n} w_i \mathbf{1}_{\{Y_i \leq y\}} + \left(1 - \frac{1}{\nu(n)} \sum_{i=1}^{n} w_i\right) \\ &= 1 - \frac{1}{\nu(n)} \sum_{i=1}^{n} w_i \mathbf{1}_{\{Y_i > y\}} \end{aligned} \tag{4.29}$$

is the proper choice. We denote the corresponding quantile estimators by

$$q^{l}_{n,w,\nu}(\alpha) = F^{l\leftarrow}_{n,w,\nu}(\alpha), \qquad q^{r}_{n,w,\nu}(\alpha) = F^{r\leftarrow}_{n,w,\nu}(\alpha). \tag{4.30}$$

The functions (4.26), (4.28) and (4.29) are no longer genuine empirical distribution functions because conditions $\lim_{x\to-\infty} F(x) = 0$ and $\lim_{x\to\infty} F(x) = 1$ may be violated. However, we still have

$$\lim_{y\to-\infty} F_{n,w,\nu}(y) = 0, \qquad \lim_{y\to-\infty} F^{l}_{n,w,\nu}(y) = 0, \qquad \lim_{y\to\infty} F^{r}_{n,w,\nu}(y) = 1.$$

For studying the convergence of the weighted quantile estimators, we assume that the sequence $(\theta_n)_{n\geq 0}$ is generated by any adaptive algorithm as described in Section 4.2 which converges to some limit value $\theta^*$. We would like to point out that it is not required that $\theta^*$ is the solution of a variance minimization problem such as given by (3.2). Later, we will propose a specific tuning algorithm and state verifiable conditions that guarantee its convergence.

ASSUMPTION 4.1. $(\mathcal{K}_j)_{j\in\mathbb{N}}$ is a compact exhaustion of the parameter space as in (4.14). The sequence $(\theta_n)_{n\geq 0}$ satisfies

$$\theta_n \to \theta^* \in \boldsymbol{\Theta} \qquad \text{almost surely.} \tag{4.31}$$

For any $\rho \in (0,1)$, there exists a constant $C(\rho)$ such that

$$\mathbb{P}\Big(\sup_{n\geq 1} \kappa_n \geq j\Big) \leq C(\rho)\rho^j, \tag{4.32}$$

where $\kappa_n$ is the counter of the active truncation set of $(\theta_n)_{n\geq 0}$ defined in such a way that (4.17) holds. For some $p^* > 1$, there exists $W \in \mathcal{L}_{p^*}(\theta_0)$ such that for any compact set $\mathcal{K} \subset \boldsymbol{\Theta}$,

$$\mathbf{1}_{\{\theta\in\mathcal{K}\}} w_x(\theta) \leq C_{p^*}(\mathrm{diam}(\mathcal{K}))W(x), \tag{4.33}$$

where $C_{p^*}(\mathrm{diam}(\mathcal{K}))$ is a constant only depending on $p^*$ and the diameter of $\mathcal{K}$. The compact covering (4.14) is selected such that

$$C_{p^*}(\mathrm{diam}(\mathcal{K}_j)) \leq e^{k_{p^*} + m_{p^*} j} \tag{4.34}$$

for some positive constants $k_{p^*}$, $m_{p^*}$.

Because of (4.32), the number of truncations remains almost surely finite and every path of $\theta_n$ remains in a compact subset of $\boldsymbol{\Theta}$. However, this does not imply that there exists a compact set $K^*$ such that $\theta_n \in K^*$ almost surely for all $n$.[5] Condition (4.33) guarantees the continuity of moments as a function of the parameters $\theta$. Condition (4.34) provides a growth restriction on the compact exhaustion (4.14).

We first address convergence when quantiles are unique but without imposing differentiability of the distribution function at the quantiles.

THEOREM 4.1. *Assume that the distribution function $F(y) = \mathbb{P}(Y \leq y)$ is strictly increasing at $q_\alpha$. Under Assumption 4.1,*

$$q_{n,w,\nu}(\alpha) \to q_\alpha \qquad \textit{almost surely } (n \to \infty),$$

*for the normalization function $\nu(n) \equiv 1$, and*

$$q^r_{n,w,\nu}(\alpha) \to q_\alpha, q^l_{n,w,\nu}(\alpha) \to q_\alpha \qquad \textit{almost surely } (n \to \infty),$$

*for $\nu(n) = n$.*

If the quantiles are not unique, a proper choice of the normalization function $\nu(n)$ eliminates the oscillatory behavior and leads to consistent estimators. For notational convenience, let

$$v = \sigma^2_{\mathbf{1}}(\theta^*) \quad \text{and} \quad v_\alpha = \sigma^2_{\mathbf{1}_{(-\infty,q_\alpha]}\circ\Psi}(\theta^*) = \sigma^2_{\mathbf{1}_{(q_\alpha,\infty)}\circ\Psi}(\theta^*). \tag{4.35}$$

---

[5] We would like to thank the anonymous referee for pointing this out to us.

THEOREM 4.2. *Suppose the conditions in Assumption 4.1 are satisfied. If there exists* $\eta > 0$, $k > 0$, *and* $0 < \gamma < \frac{1}{2}$ *such that*

$$\frac{n - kn^{1/2+\gamma}}{n - (1+\eta)\sqrt{2nv \log\log(nv)}} \le \nu(n) \le \frac{n - (1+\eta)/\alpha\sqrt{2nv_\alpha \log\log(nv_\alpha)}}{n + (1+\eta)\sqrt{2nv \log\log(nv)}}, \tag{4.36}$$

*then*

$$q_{n,w,\nu}(\alpha) \to q_\alpha \qquad \textit{almost surely } (n \to \infty). \tag{4.37}$$

*If there exists* $\eta > 0$, $k > 0$, *and* $0 < \gamma < \frac{1}{2}$ *such that*

$$n + \frac{1+\eta}{1-\alpha}\sqrt{2nv_\alpha \log\log(nv_\alpha)} \le \nu(n) \le n + kn^{1/2+\gamma}, \tag{4.38}$$

*then*

$$q^r_{n,w,\nu}(\alpha) \to q_\alpha \qquad \textit{almost surely } (n \to \infty). \tag{4.39}$$

*If there exist* $\eta > 0$, $k > 0$, *and* $0 < \gamma < \frac{1}{2}$ *such that*

$$n - kn^{1/2+\gamma} \le \nu(n) \le n - \frac{1+\eta}{\alpha}\sqrt{2nv_\alpha \log\log(nv_\alpha)}, \tag{4.40}$$

*then*

$$q^l_{n,w,\nu}(\alpha) \to q_\alpha \qquad \textit{almost surely } (n \to \infty). \tag{4.41}$$

The proofs of Theorems 4.1 and 4.2 are given in Section 6. They rely on a law of iterated logarithm for martingales which we present in a later section (Section 4.6).

The normalization functions used in Theorem 4.2 are difficult to implement, because $v_\alpha$ depends on the unknown quantile $q_\alpha$ and the unknown limit parameter $\theta^*$. In this regard, the following corollary is helpful.

COROLLARY 4.1. *If* $\theta^* = \arg\min_\theta \sigma^2_{\mathbf{1}_{(q_\alpha,\infty)}\circ\Psi}(\theta)$, *then*

$$\begin{aligned} v_\alpha &= \sigma^2_{\mathbf{1}_{(q_\alpha,\infty)}\circ\Psi}(\theta^*) \\ &\le \sigma^2_{\mathbf{1}_{(q_\alpha,\infty)}\circ\Psi}(\theta_0) = \mathbb{P}_{\theta_0}(\Psi(X) > q_\alpha) - \mathbb{P}_{\theta_0}(\Psi(X) > q_\alpha)^2 \le \tfrac{1}{4}. \end{aligned}$$

Therefore, the conclusions of Theorem 4.2 hold, if $v_\alpha$ in conditions (4.36), (4.38), and (4.40) is replaced by $\frac{1}{4}$.

To compare Theorem 4.2 with Theorem 4 of Feldman and Tucker in [12], we state here a refined version of their result.

THEOREM 4.3. *Let $Y_{1,n},\dots,Y_{n,n}$ be the order statistics of i.i.d. samples $Y_1,\dots,Y_n$ of a random variable $Y$. Let*

$$w_\alpha = \mathbb{P}(Y \le q_\alpha(Y)) - \mathbb{P}(Y \le q_\alpha(Y))^2. \tag{4.42}$$

*If the normalization function $\nu(n) \in \mathbb{N}$ satisfies*

$$(1+k)\sqrt{2w_\alpha n \log\log n} \le \lfloor n\alpha \rfloor - \nu(n) \le K n^{1/2+\gamma} \tag{4.43}$$

*with $\gamma, k, K$ positive constants, then $Y_{\nu(n),n} \to q_\alpha(Y)$ almost surely.*

We omit the proof, as it is similar to the proof of Theorem 4.2. Condition (4.43) is now expressed in a way that allows a direct comparison with (4.38). We recall the original condition in Theorem 4 of Feldman and Tucker,

$$(1+k)\sqrt{2n \log\log n} \le \lfloor n\alpha \rfloor - \nu(n) \le K n^{1/2+\gamma}, \tag{4.44}$$

which apparently does not depend on the variance of the tail probabilities. However, because $w_\alpha \le \frac{1}{4}$ we see that (4.43) is indeed a weaker assumption than (4.44) used in [12].

4.4. *Sequential parameter tuning for quantile estimation.* We still have to provide a strategy to determine the limit parameter $\theta^*$ and the construction of an approximation sequence $\theta_n$ converging to $\theta^*$ almost surely. For the estimation of the expected value $\mathbb{E}[f(X)]$, Theorem 3.1 suggests that the optimal parameter $\theta^*$ is the variance minimizer $\theta^* = \arg\min_\theta \sigma_f^2(\theta)$ which can be estimated by a stochastic approximation algorithm. However, for quantile estimation, the choice of an optimal parameter $\theta^*$ is less obvious. If the distribution function $F$ of the random variable $Y = \Psi(X)$ is differentiable in a neighborhood of $q_\alpha$, the functional delta-method applied to the empirical process (see, e.g., Corollary 21.5 of [33]) suggests to minimize the variance of the weighted tail event $w_X(\theta)\mathbf{1}_{\{\Psi(X)>q_\alpha\}}$ such that

$$\theta^* = \arg\min_\theta m_{\mathbf{1}_{(q_\alpha,\infty)}\circ\Psi}(\theta). \tag{4.45}$$

Instead of arguing with the delta-method as above, we can also use Theorem 4.2 to motivate the choice (4.45) even in the most general situation, in which quantiles may not be unique. For instance, let us consider the quantile estimator $q^r_{n,w,\nu}$. The bounds for $\nu(n)$ in (4.38) lead to a bias for $q^r_{n,w,\nu}$. To minimize this bias, we must ensure that $\nu(n)$ is as close as possible to $n$ while, at the same time, satisfying condition (4.38). This means that we must select $v_\alpha$ such that the term $\sqrt{2nv_\alpha \log\log(nv_\alpha)}$ becomes as small as possible. From the definition of $v_\alpha$ in (4.35), we see that the parameter $\theta^*$ satisfying (4.45) provides the smallest value for $v_\alpha$. The same arguments hold for $q^l_{n,w,\nu}$.

For the estimator $q_{n,w,\nu}(\alpha)$ defined in (4.27), we must keep $\nu(n)$ as close as possible to 1 in order to minimize the bias. From condition (4.36), we see that we must not only minimize $\sqrt{2nv_\alpha \log\log(nv_\alpha)}$, but also $\sqrt{2nv\log\log(nv)}$. Hence, for $q_{n,w,\nu}(\alpha)$ we have to choose $\theta$ to make both the variance of the weighted tail event $w_X(\theta)\mathbf{1}_{\{\Psi(X)>q_\alpha\}}$ and the variance of the weights $w_X(\theta)$ as small as possible.

Unfortunately (4.45) is not constructive either because the quantile $q_\alpha$ is not yet known and must be replaced by a suitable estimator. Suppose now that we could find a rough estimate $\hat{q}_\alpha$ for the quantile $q_\alpha$; then the scheme (4.5) based on the stochastic gradient,

$$H_{\hat{q}_\alpha}(X_{n+1},\theta_n), \qquad X_{n+1}\sim\varphi_{\theta_n}\,d\lambda \tag{4.46}$$

with

$$H_q(x,\theta)=-\mathbf{1}_{\{\Psi(x)>q\}}w_x(\theta)^2\nabla l_x(\theta) \tag{4.47}$$

could be used to generate a sequence $(\theta_n)_{n\geq 0}$ approximating the solution $\theta^*$ for the first order condition $\nabla m_{\mathbf{1}_{(q_\alpha,\infty)}\circ\Psi}(\theta^*)=0$. However, if $q_\alpha$ is an extreme quantile, the simulated values for the stochastic gradient (4.46) would be mostly zero for parameter values $\theta_n$ close to the starting value $\theta_0$. Even worse, if the simulation produces a nonvanishing stochastic gradient, it is generally very inaccurate and could drive the parameter values to a wrong region of the parameter space. As a consequence, the convergence rate of the algorithm is very poor. It freezes at an early stage and one might be tempted to use a sufficiently large step size. However, in practical applications, compensating an erratic stochastic gradient with a large step size is not a solution, as it increases the risk that the algorithm fails to converge.

A simple and practically very efficient approach is to gradually bridge from a moderate tail event to an extreme tail event during the simulation. More precisely, let

$$M_{q_1,q_2}(\theta)=b(n)m_{\mathbf{1}_{(q_1,\infty)}\circ\Psi}(\theta)+(1-b(n))m_{\mathbf{1}_{(q_2,\infty)}\circ\Psi}(\theta) \tag{4.48}$$

with $b(n)$ weighting functions depending on the sample index $n$. The values $q_i$ are selected such that $q_\alpha\in[q_1,q_2]$. We choose $q_1$ such that $\{\Psi(X)>q_1\}$ is a moderate tail event. Hence, the corresponding stochastic gradient $H_{q_1}(X_{n+1},\theta_n)$ can be estimated with sufficient accuracy for $\theta_n$ in a neighborhood of $\theta_0$. The value $q_2$ is selected in the range of $q_\alpha$ or even larger. A preliminary simulation or some initial samples can be used to obtain a crude estimate for $q_\alpha$, including a confidence interval. The function $b(n)$ is assumed to converge to zero as $n\to\infty$. A suitable choice would be, for example, $b(n)=1/\log(n+1)$ which decays sufficiently slow such that the component (4.46) of the stochastic gradient from $q_1$ drives $\theta_n$ towards a solution for the extreme tail event. Stochastic approximation with adaptive truncation can then be used to generate a sequence of parameters $\theta_n$ converging to

$$\theta^* = \arg\min_{\theta} M_{q_1,q_2}(\theta) \tag{4.49}$$

as we will see below.[6]

4.5. *Verifiable convergence criteria.* Each of the above criterion is based on a stochastic vector field generated by a map $H(x,\theta)$. For instance, in case of (4.48), we have

$$H(x,\theta) = b(n)H_{q_1}(x,\theta) + (1-b(n))H_{q_2}(x,\theta). \tag{4.50}$$

We provide verifiable conditions on $H(x,\theta)$, its mean field, and the sequences $\boldsymbol{\gamma}$, $\boldsymbol{\epsilon}$, which imply the convergence of Algorithm 4.1 for state independent transition kernels. To this end, we introduce for any compact set $\mathcal{K} \subset \boldsymbol{\Theta}$ the partial sum

$$S_{l,n}(\boldsymbol{\gamma},\boldsymbol{\epsilon},\mathcal{K}) = \mathbf{1}_{\{\sigma(\mathcal{K},\boldsymbol{\epsilon})\geq n\}} \sum_{k=l}^{n} \gamma_k (H(X_k,\theta_{k-1}) - h(\theta_{k-1})), \qquad 1 \leq l \leq n, \tag{4.51}$$

where $\sigma(\mathcal{K},\boldsymbol{\epsilon}) = \sigma(\mathcal{K}) \wedge \nu(\boldsymbol{\epsilon})$ and $\sigma(\mathcal{K})$ and $\nu(\boldsymbol{\epsilon})$ are the stopping times

$$\sigma(\mathcal{K}) = \inf\{k \geq 1 \mid \theta_k \notin \mathcal{K}\}, \tag{4.52}$$

$$\nu(\boldsymbol{\epsilon}) = \inf\{k \geq 1 \mid |\theta_k - \theta_{k-1}| \geq \epsilon_k\}. \tag{4.53}$$

If $\mathbf{a} = (a_l)_{l\in\mathbb{N}}$ is a sequence, we write

$$\mathbf{a}^{\leftarrow k} = (a_{l+k})_{l\in\mathbb{N}}$$

for the sequence shifted by the offset $k$.

ASSUMPTION 4.2. The parameter set $\boldsymbol{\Theta}$ is an open subset of $\mathbb{R}^d$. For some $p > 1$, there exists a function $W \in \mathcal{L}_p(\theta_0)$ such that for every compact set $\mathcal{K} \subset \boldsymbol{\Theta}$,

$$\sup_{x\in\mathcal{X}} \sup_{\theta\in\mathcal{K}} \frac{\|H(x,\theta)\|^p}{w_x(\theta)W(x)^p} \leq C_{\mathcal{K}} < \infty \tag{4.54}$$

[6] Yet another approach is to sequentially update an estimator $\hat{q}_\alpha$ for the quantile along the simulation as well to improve the upper value $q_2$ in (4.49), leading to a coupled stochastic approximation scheme for the parameters $(\theta_n, \hat{q}_n)$. A sequential quantile estimator has been proposed in [32] (see also [31]). Because the quantile estimator does interfere with the update scheme for the tuning parameter $\theta_n$, the convergence of the joint parameter set is more subtle.

with $C_{\mathcal{K}}$ a constant only depending on $\mathcal{K}$. The mean field

$$h(\theta) = \mathbb{E}_\theta[H(X,\theta)] \tag{4.55}$$

is continuous and there exists a $C^1$ Lyapunov function $w:\boldsymbol{\Theta} \to [0,\infty)$ satisfying the following conditions:

(i) There exists $0 < M_0 < \infty$ such that

$$\mathcal{L} \equiv \{\theta \in \boldsymbol{\Theta} \mid \langle h(\theta), \nabla w(\theta)\rangle = 0\} \subset \{\theta \in \boldsymbol{\Theta} \mid w(\theta) < M_0\}.$$

(ii) For $M > 0$, let $\mathcal{W}_M = \{\theta \in \boldsymbol{\Theta} \mid w(\theta) \le M\}$. There exists $M_1 \in (M_0, \infty]$ such that $\mathcal{W}_{M_1}$ is a compact subset of $\boldsymbol{\Theta}$.

(iii) For any $\theta \in \boldsymbol{\Theta} \setminus \mathcal{L}$ it holds that $\langle h(\theta), \nabla w(\theta)\rangle < 0$.

(iv) The closure of $w(\mathcal{L})$ has empty interior.

The sequences $\boldsymbol{\gamma}$, $\boldsymbol{\epsilon}$ are nonincreasing, positive, and satisfy $\epsilon_n \to 0$,

$$\sum_{n=0}^{\infty} \gamma_n = \infty, \qquad \sum_{n=0}^{\infty} \left( \gamma_n^2 + \left(\frac{\gamma_n}{\epsilon_n}\right)^p \right) < \infty. \tag{4.56}$$

The existence of a Lyapunov function in (i) simplifies, if $h = \nabla m$ is a gradient field of a continuously differentiable function $m$. In this case, we can choose $w = m$. The next result is along the lines of Proposition 5.2 in [1]. Its proof is similar to the proof of [1], Proposition 5.2, but less involved because we consider only state-independent transition probabilities. Therefore, we do not need to consider the existence and regularity of the solution of the Poisson equation. The convergence of the algorithm is then a consequence of [1], Theorem 5.5.

Proposition 4.1. *Let*

$$\begin{aligned} A(\delta, M, \boldsymbol{\gamma}, \boldsymbol{\epsilon}) = \sup_{(x,\theta)\in\mathcal{X}_0\times\mathcal{K}_0} \Big\{ &\mathbb{P}^{\boldsymbol{\gamma}}_{\Phi(x,\theta)}\Big(\sup_{n\ge 1} \|S_{1,n}(\boldsymbol{\gamma}, \boldsymbol{\epsilon}, \mathcal{W}_M)\| > \delta\Big) \\ &+ \mathbb{P}^{\boldsymbol{\gamma}}_{\Phi(x,\theta)}(\nu(\boldsymbol{\epsilon}) < \mathcal{W}_M)\Big\}. \end{aligned} \tag{4.57}$$

*If $\mathcal{K}_0 \subset \mathcal{W}_{M_0}$, then for every $M \in [M_0, M_1)$ there exist $n_0$, $\delta_0 > 0$, and a constant $C > 0$ such that for all $j > n_0$,*

$$\mathbb{P}\Big(\sup_{k\ge 1} \kappa_k \ge j\Big) \le C\Big(\sup_{k \ge n} A(\delta_0, M, \boldsymbol{\gamma}^{\leftarrow k}, \boldsymbol{\epsilon}^{\leftarrow k})\Big)^j. \tag{4.58}$$

*Under Assumption 4.2, we have for every $M \in [M_0, M_1)$ and $\delta > 0$,*

$$\lim_{k\to\infty} A(\delta, M, \boldsymbol{\gamma}^{\leftarrow k}, \boldsymbol{\epsilon}^{\leftarrow k}) = 0. \tag{4.59}$$

*In particular, the key requirements, (4.31) and (4.32), of Assumption 4.1 are satisfied.*

To completely specify the stochastic approximation algorithm, we first have to make some selections for the initial parameter $\theta_0$. Because our target criterion puts more emphasis on a moderate tail event at the beginning of the simulation, it is sensible to start with the reference density. Alternatively, we can start with a large deviation approximation.

The performance of a stochastic approximation algorithm usually depends strongly on an appropriate selection of the step size sequence. However, with the bridging strategy in (4.48), our algorithm is considerably less sensitive to the choice of the step size parameters. Since the sequence of step size parameters $\gamma_n$ must satisfy condition (4.56), we simply set

$$\gamma_n = \frac{a}{n+1} \tag{4.60}$$

and select $\epsilon_n$ accordingly to satisfy the second condition in (4.56). The parameter $a$ serves as a tuning parameter. A practical approach is to follow a greedy strategy which starts with a large value for $a$ and reduces it after each re-initialization. Alternatively we can determine $a$ by some step size selection criteria based, for example, on an approximation of the Hessian of the target criterion.

The algorithm can be further robustified by Polyak's averaging principle. The idea is to use a large step size $\gamma_n$ of the order $n^{-2/3}$ which converges much slower to zero than $n^{-1}$ but is still fast enough to ensure convergence. The larger step size prevents the algorithm from freezing at an early stage of the algorithm far off the local minimum. Polyak and Juditsky show in [27] that the averaged parameters converge at an optimal rate.

4.6. *Law of iterated logarithm for martingales.* Before we discuss an application for our adaptive quantile estimator, we present the law of iterated logarithm for the sequence of martingale differences $w_{X_n}(\theta_{n-1})f(X_n) - \mathbb{E}[f(X)]$ which we require as an ingredient for the proofs of Theorems 4.1 and 4.2. We state the main result below and present the proof in Section 6. We use the following notation. If $M_n$ is a square integrable martingale adapted to a filtration $(\mathcal{F}_n)_{n\geq 0}$ with $\Delta M_i = M_i - M_{i-1}$, then we denote the predictable quadratic variation by

$$\langle M\rangle_0 = 0, \qquad \langle M\rangle_n = \sum_{i=1}^{n} \mathbb{E}[\Delta M_i^2 \mid \mathcal{F}_{i-1}], \qquad n \geq 1, \tag{4.61}$$

the total quadratic variation by

$$[M]_0 = 0, \qquad [M]_n = \sum_{i=1}^{n} \Delta M_i^2, \qquad n \geq 1, \tag{4.62}$$

and by $s_n^2 = \sum_{i=1}^n \mathbb{E}[\Delta M_i^2]$ the total variance.

THEOREM 4.4. *Suppose the conditions in Assumption 4.1 are satisfied, and let $f:\mathcal{X}\to\mathbb{R}$ be a measurable function in $\mathcal{L}_p(\theta_0)$. Assume that*

$$\frac{p(p^*+1)}{p+p^*}>4, \tag{4.63}$$

*where $p^*$ is from condition (4.33). Let*

$$w_x:\theta\mapsto w_x(\theta) \tag{4.64}$$

*be continuous in $\theta^*$ for almost all $x\in\mathcal{X}$. Define*

$$\xi_n=w_{X_n}(\theta_{n-1})f(X_n)-\mathbb{E}[f(X)]. \tag{4.65}$$

*Then $M_n=\sum_{i=1}^n\xi_i$ is a square integrable martingale and*

$$\lim_{n\to\infty}\frac{[M]_n}{\langle M\rangle_n}=1, \tag{4.66}$$

$$\lim_{n\to\infty}\frac{s_n^2}{n}=(m_f(\theta^*)-\mathbb{E}[f(X)]^2)=\sigma_f^2(\theta^*). \tag{4.67}$$

*Moreover, if we let $\phi(t)=\sqrt{2t\log\log(t)}$, then*

$$\limsup_{n\to\infty}\phi(W_n)^{-1}M_n=+1, \tag{4.68}$$

$$\liminf_{n\to\infty}\phi(W_n)^{-1}M_n=-1, \tag{4.69}$$

*almost surely, where the weighting sequence $W_n$ is given by either $W_n=[M]_n$, $W_n=\langle M\rangle_n$, or $W_n=s_n^2$.*

**5. Applications.** We next provide an explicit example for our adaptive quantile estimator and compare it to crude Monte Carlo simulation. We borrow our application from the financial industry, more precisely from portfolio credit risk. The so-called Value at Risk (VaR) is by far the most widely adopted measure of risk and represents the maximum level of losses that can be exceeded only with a small probability. This quantile-based risk measure is of particular importance to market participants and supervisors. For credit risk, supervisors require banks to calculate the credit VaR as the 99.9% quantile of the loss distribution.

5.1. *Importance sampling for portfolio credit risk.* The aim of portfolio credit risk analysis is to provide a distribution of future credit losses for a portfolio of obligers based on historically observed losses and possibly combined with market views. In a simplified setting, the outstanding credit amount for each obligor $i=1,\dots,m$ is aggregated to a net credit exposure

$c_i$. Defaults are tracked over a single period. At the end of the period, the portfolio loss is

$$\mathcal{L} = \sum_{i=1}^{m} c_i Y_i, \tag{5.1}$$

where $Y_i \sim \mathrm{Ber}(p_i)$ are the default indicators. For portfolios of illiquid commercial loans or corporate credits, the exposures $c_i$ are generally assumed to be constant which gives rise to a discrete loss distribution. The quantiles are nondifferentiable and not unique. Hence, to construct an adaptive importance sampling algorithm, we can rely on the results of the previous sections.

For our application, we start from a Gaussian copula framework (see, e.g., [11]), in which the default indicators are modeled as

$$Y_i = \mathbf{1}_{\{A_i \in (-\infty, \theta_i]\}}. \tag{5.2}$$

The credit quality variable $A_i$ is given by

$$A_i = \sqrt{1 - v_{s(i)}^2} X_{s(i)} + v_{s(i)} \epsilon_i, \qquad i = 1, \ldots, m, \tag{5.3}$$

for some classification function $s : \{1, \ldots, m\} \to \{1, \ldots, k\}$. Usually in credit risk management, the $m$ obligors are classified into $k$ industry sectors. The default thresholds $\theta_i$ are calibrated to match the obligors' default probabilities. The common factors $X = (X_s)_{s=1,\ldots,k} \sim N(0, \Sigma)$ are multivariate Gaussian. The idiosyncratic part $\epsilon = (\epsilon_i)_{i=1,\ldots,m} \sim N(0, \mathbf{1}_m)$ is independent from $X$. We restrict ourselves to the adjustment of the mean of the common factors $X$ and keep the covariance structure $\Sigma$ fixed. We note that importance sampling on the common factors can also be combined seamlessly with importance sampling on the idiosyncratic variables $\epsilon = (\epsilon_i)_{i=1,\ldots,m}$.[7]

Given the above setup, we are in the setting of Section 4.3 with $Y = \mathcal{L} = \Psi(X)$ and $\Psi : \mathbb{R}^k \to \mathbb{R}$ given as[8]

$$\Psi(x) = \sum_{i=1}^{m} c_i \mathbf{1}_{\{\sqrt{1 - v_{s(i)}^2} x_{s(i)} + v_{s(i)} \epsilon_i \le \theta_i\}}. \tag{5.4}$$

For the implementation of the adaptive importance sampling scheme, we use the criterion (4.48) and determine the values for $q_1$, $q_2$ as described in Section 4.4, that is, we start with a moderate $q_1$ and choose $q_2$ by an educated guess in the region of interest.

---

[7] For instance, [13] and [26] apply an exponential twist to the conditional default indicators $Y_i \mid X \sim \mathrm{Ber}(p_i(X))$ where $p_i(x) = \Phi((\theta_i - \sqrt{1 - v_{s(i)}^2} x_{s(i)})/v_{s(i)})$, is the conditional default probability.

[8] For notational convenience, we drop the dependency of $\Psi$ on $\epsilon$ as it is not affected by the importance sampling scheme.

5.2. *Verifying convergence criteria for Gaussian distributions.* For our credit risk application, assume a fixed covariance structure, and endow the Gaussian distributions (4.9) with the Fisher information metric $g_\Sigma$ in (4.11). Before we can proceed, we need to make sure that Assumptions 4.1 and 4.2 hold in our setup. The noncompactness of the parameter space and exponentially unbounded likelihood ratios call for adaptive truncation and make it a challenging test case, even though the Gaussian distributions have many special analytical properties. From the expression

$$w_x(\theta) = \exp(-g_\Sigma(x,\theta) + \tfrac{1}{2} g_\Sigma(\theta,\theta))$$

for the likelihood ratio, it follows that

$$w_x(\theta) \le \exp\left(\frac{h+2}{4} g_\Sigma(\theta,\theta)\right) \exp\left(\frac{1}{h} g_\Sigma(x,x)\right) \qquad \forall h \ge 1. \tag{5.5}$$

The verification of Assumptions 4.1 and 4.2 for the Gaussian distributions is now a straightforward consequence of (5.5) and Hölder's inequality.

LEMMA 5.1. *If $\|f\|_{\theta_0,h} < \infty$ for some $h > 2$, we can exchange differentiation and integration to obtain $\nabla m_f(\theta) = \mathbb{E}_\theta[H(X,\theta)]$ with*

$$H(x,\theta) = (\theta - x) f(x)^2 w_x(\theta)^2.$$

*The Hessian with respect to the Fisher information metric $g_\Sigma$, given by*

$$\nabla^2 m_f(\theta) = \mathbb{E}_\theta[(\mathrm{id}_k + (\theta - X)(\theta - X)^\top) f(X)^2 w_X(\theta)^2], \tag{5.6}$$

*is positive definite. If $\mathbb{P}(f(X) > 0) > 0$, then $m_f(\theta) \to \infty$ for $g_\Sigma(\theta,\theta) \to \infty$. In particular, there is a unique minimizer*

$$\theta^* = \arg\min_\theta m_f(\theta) \in \mathbb{R}^k. \tag{5.7}$$

*Moreover, for some $p > 1$ there exists $W \in \mathcal{L}_p(\theta_0)$ satisfying (4.33) and (4.54).*

The parametrization (4.9) works rather well if the ratio of the largest and smallest eigenvalue of $\Sigma$ is not too far away from one and the dimension of $\Sigma$ is not too large. For many practical applications, the correlation ellipsoid is very skewed. The first few principal components explain most of the variance and the last few are negligibly small. Even though the metric defined in (4.11) properly respects the covariance structure, and we use the gradient relative to this metric, we require a suitable dimension reduction. Therefore, we translate the mean in the span of the eigenvalues of the first few principal components. Let $\Sigma = U \Lambda U^T$ where $U$ is the orthogonal matrix with columns given by the eigenvectors, and $\Lambda$ is the diagonal matrix of eigenvalues. We write

$$J_l : \mathbb{R}^l \to \mathbb{R}^k \tag{5.8}$$

for the embedding of $\mathbb{R}^l$ into $\mathbb{R}^k$; that is, $J_l$ sets the last $k-l$ coordinates to zero with corresponding projection $J_l^\top : \mathbb{R}^k \to \mathbb{R}^l$. Let

$$\mathbf{\Theta}_l = \{N(UJ_l(a), \Sigma) \mid a \in \mathbb{R}^l\}. \tag{5.9}$$

The first and second order differential of the likelihood $l_x(a)$ is

$$dl_x(a) = J_l^\top \Lambda^{-1}(U^\top x - J_l a), \qquad d^2 l_x(a) = -J^\top \Lambda^{-1} J. \tag{5.10}$$

Hence, the Fisher metric on $\mathbf{\Theta}_l$ is

$$g_a(u, v) = -\mathbb{E}_a[d^2 l_X(a)(u, v)] = u^\top J^\top \Lambda^{-1} J v. \tag{5.11}$$

Because

$$g_a(\nabla l_x(a), u) = \nabla l_x(a)^\top J_l^\top \Lambda^{-1} J_l u = dl_x(a)(u) = (x^\top U - a^\top J_l^\top)\Lambda^{-1} J_l u$$

and $x^\top U \Lambda^{-1} J_l u = x^\top U J_l J_l^\top \Lambda^{-1} J_l u$, the gradient of the likelihood with respect to the metric (4.11) is

$$\nabla l_x(a) = (J_l^\top U^\top x - a). \tag{5.12}$$

We adapt Lemma 5.1 to the parametrization given in (5.9).

LEMMA 5.2. *If $\|f\|_{\theta_0,q} < \infty$ for some $q > 2$, we can exchange differentiation and integration to obtain $\nabla m_f(a) = \mathbb{E}_{\theta(a)}[H(X, a)]$ with*

$$H(x, a) = (a - J_l^\top U^\top x) f(x)^2 w_x(\theta(a))^2$$

*and $\theta(a) = UJ_l(a)$. The Hessian with respect to the Fisher information metric $g_a$, given by*

$$\nabla^2 m_f(\theta) = \mathbb{E}_\theta[(\mathrm{id}_l + (a - J_l^\top U^\top X)(a - J_l^\top U^\top X)^\top) f(X)^2 w_X(\theta(a))^2],$$

*is positive definite. If $\mathbb{P}(f(X) > 0) > 0$, then $m_f(a) \to \infty$ for $g_a(a, a) \to \infty$. In particular, there is a unique minimizer*

$$a^* = \arg\min_a m_f(a) \in \mathbb{R}^l. \tag{5.13}$$

*Moreover, for some $p > 1$ there exists $W \in \mathcal{L}_p(\theta_0)$ satisfying (4.33) and (4.54).*

5.3. *Numerical example.* We consider a set of 2000 obligors with default probabilities comparable to a typical loan portfolio. We assume that the portfolio risk is driven by $m = 14$ industry factors, but restrict our analysis using only the first two principal components which already explain 84% of total variance. In the current regulatory framework as promoted by Basel II, credit risk (as well as operational risk) needs to be calculated at the 99.9% quantile of the loss distribution. Performing a crude Monte Carlo (MC) simulation, we see that the loss (expressed in percentage numbers) at the 99.9% quantile lies somewhere around 0.2. This crude estimate allows us to make an educated guess for the parameters $q_1$ and $q_2$ required for our adaptive importance sampling (AIS) estimator. We set $q_1 = 0.1$ and $q_2 = 0.23$. Instead of using the MC estimate as a starting point, we could also first do an AIS simulation with some arbitrarily set $q_1$ and $q_2$ to find some more appropriate numbers in a second simulation. Our numerical experiments indicated that the algorithm is not very sensitive to these approximate choices. Indeed, we just have to guarantee that we choose $q_1$ small so that the initial step sizes are large enough. To clarify this point with an example, we find that we get almost identical results for $q_1 = 0.01$. More precisely, with a fixed seed for the random number generator we get an estimate for 99.9%-quantile of 0.2271 with $q_1 = 0.1$ and 0.2276 with $q_1 = 0.01$.

Based on a sample of 10,000 draws, Figure 1 shows the convergence of the mean shifts for our AIS algorithm. The solid line represents the path for the step size of order $n^{-2/3}$. The dashed line represents the averaged values based on Polyak's averaging principle. We observe that the shift in the first principal component, which explains 75% of total variance, is substantial. The shift in the second component, which explains an additional 9%, is only very small.

Figure 2 plots the cumulative distribution function for the right tail of the distribution. In contrast to standard MC simulation, our AIS algorithm provides a very smooth distribution function. Therefore, we can expect a considerable reduction for the variance of our quantile estimators. To substantiate this conjecture, we additionally perform 1000 independent quantile estimations. In Table 1, we report the results for the standard Monte Carlo simulations to calculate $F_{n,1}^{\leftarrow}(\alpha)$ and for our AIS algorithm using the quantile definition in (4.30) which is based on the weighted empirical distribution $F_{n,w,\nu}^{r}(y)$. The first column shows the different loss levels at which we simulate the quantiles. The next two columns report the mean values for the estimations $F_{n,1}^{\leftarrow}(\alpha)$ and $F_{n,w,\nu}^{r,\leftarrow}(\alpha)$, respectively. The final row reports the variance ratio defined as the variance from the MC simulation divided by the variance of the AIS estimator. When we compare the variances of the estimators, we observe that for the region of interest, that is, around the 99.9% quantile, our AIS estimator outperforms the result from the MC simulation by a factor of around 20. This number increases further to more than 112 when we look at the 99.99% quantile.

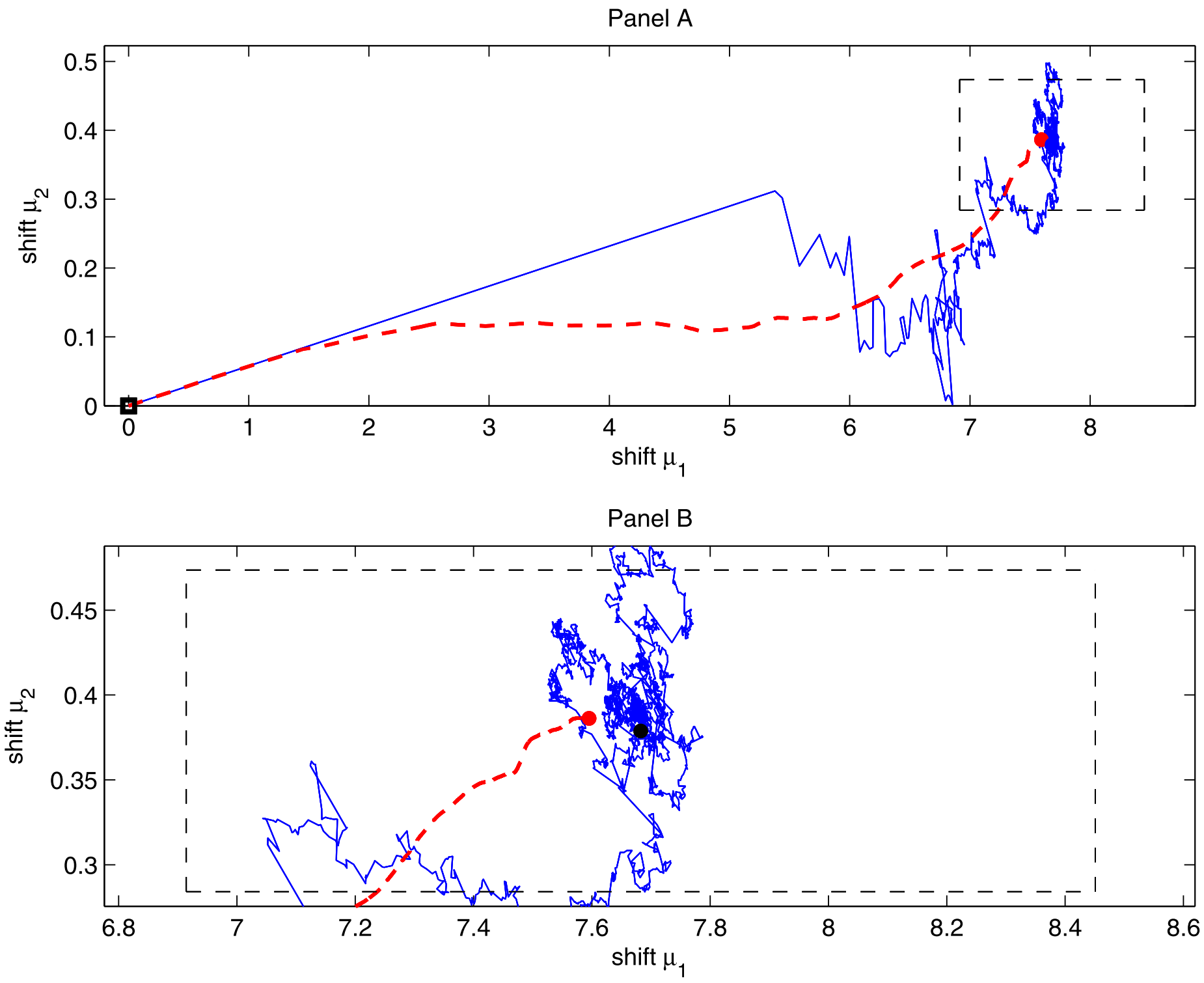

Fig. 1. *Parameter convergence from Polyak's averaging principle. The graph illustrates the convergence of the means of the first two principal components. The solid line represents the convergence for step size of order* $n^{-2/3}$. *The dashed line represents the averaged values. Panel* B *zooms in the rectangle area marked in panel* A.

**6. Proofs.** For the proofs for Theorems 4.1, 4.2 and 4.4 we start with collecting the basic properties of the generalized inverse of an increasing function.[9]

Lemma 6.1. *Let* $F$ *be a right continuous increasing function. Then, the generalized inverse*

$$F^{\leftarrow}(\alpha) = \inf\{x \mid F(x) \geq \alpha\} \tag{6.1}$$

*is increasing and left continuous, and we have*

$$F(x) \geq \alpha \quad \Leftrightarrow \quad F^{\leftarrow}(\alpha) \leq x;$$
$$F(x) < \alpha \quad \Leftrightarrow \quad F^{\leftarrow}(\alpha) > x;$$
$$F(x_1) < \alpha \leq F(x_2) \quad \Leftrightarrow \quad x_1 < F^{\leftarrow}(\alpha) \leq x_2;$$

[9]See for instance [29], Section 0.2.

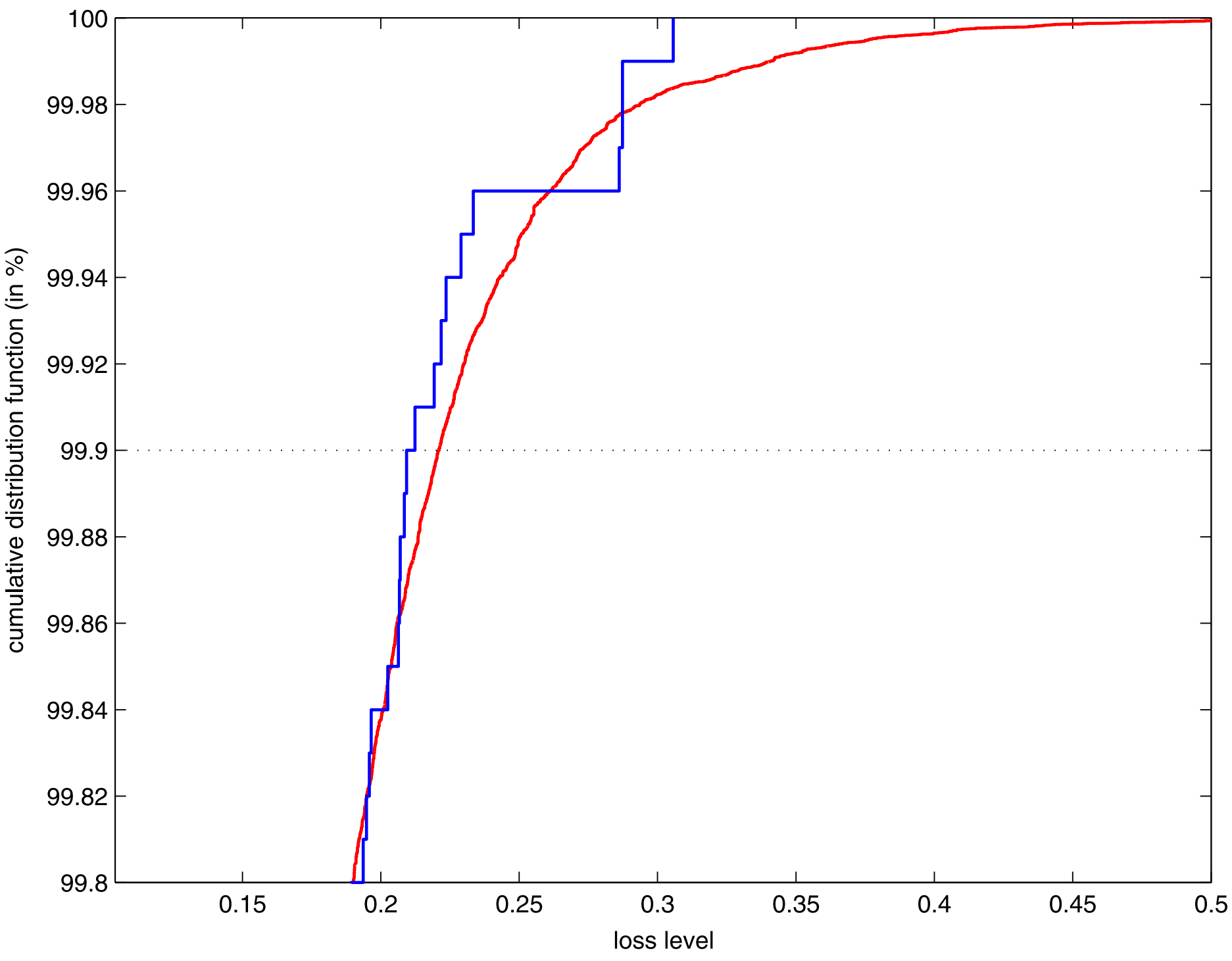

FIG. 2. *Cumulative distribution function of Monte Carlo and AIS simulation. The graph plots the cumulative distribution function using a Monte Carlo simulation and the AIS algorithm based on 10,000 samples. Our area of interest, the 99.9% quantile, is marked with a dotted line.*

$$F(F^{\leftarrow}(\alpha)) \geq \alpha, \text{ with equality for } F \text{ continuous};$$
$$F^{\leftarrow}(F(x)) \leq x, \text{ with equality for } F^{\leftarrow} \text{ increasing};$$
$$F \text{ continuous} \quad \Leftrightarrow \quad F^{\leftarrow} \text{ increasing};$$
$$F \text{ increasing} \quad \Leftrightarrow \quad F^{\leftarrow} \text{ continuous}.$$

6.1. *Proof of Theorem 4.4 (iterated law of logarithm).* Let $\kappa_n$ denote the active truncation set for $\theta_n$ and

$$\kappa_\infty = \lim_{n\to\infty} \kappa_n, \tag{6.2}$$

which exists because of assumption (4.32). We have

$$\theta_i \in \mathcal{K}_{\kappa_n} \qquad \forall i \leq n \text{ with } \nu_i \neq 0.$$

To get rid of the condition $\nu_i \neq 0$, we decompose $M_n$ into

$$M_n = \sum_{i=1}^{n} \xi_i \mathbf{1}_{\{\nu_i \neq 0\}} + \sum_{i=1}^{n} \xi_i \mathbf{1}_{\{\nu_i = 0\}}. \tag{6.3}$$

TABLE 1
*Quantile estimates for standard MC and AIS simulation at different loss levels based on 1000 simulations. Quantiles are expressed in percentage numbers, and the variance ratio is defined as the variance of $F_{n,1}^{\leftarrow}(\alpha)$ divided by the variance of $F_{n,w,\nu}^{r;\leftarrow}(\alpha)$*

| | Mean (in %) | | |
|---|---|---|---|
| **Loss level $\lambda$** | $F_{n,w,\nu}^{r;\leftarrow}(\alpha)$ | $F_{n,1}^{\leftarrow}(\alpha)$ | **Variance ratio** |
| 0.1000 | 98.2449 | 98.2217 | 1.7988 |
| 0.1100 | 98.6193 | 98.6011 | 2.0759 |
| 0.1200 | 98.9099 | 98.8943 | 2.6251 |
| 0.1300 | 99.1366 | 99.1238 | 3.0764 |
| 0.1400 | 99.3162 | 99.3053 | 3.7290 |
| 0.1500 | 99.4577 | 99.4488 | 4.6726 |
| 0.1600 | 99.5701 | 99.5638 | 5.9360 |
| 0.1700 | 99.6596 | 99.6559 | 6.0449 |
| 0.1800 | 99.7309 | 99.7287 | 7.3198 |
| 0.1900 | 99.7877 | 99.7851 | 9.7769 |
| 0.2000 | 99.8322 | 99.8305 | 12.0368 |
| 0.2100 | 99.8677 | 99.8663 | 14.3319 |
| 0.2200 | 99.8955 | 99.8946 | 17.5088 |
| 0.2300 | 99.9173 | 99.9167 | 21.1197 |
| 0.2400 | 99.9344 | 99.9340 | 25.3890 |
| 0.2500 | 99.9477 | 99.9474 | 29.9903 |
| 0.2600 | 99.9581 | 99.9578 | 35.6695 |
| 0.2700 | 99.9662 | 99.9658 | 40.5849 |
| 0.2800 | 99.9726 | 99.9724 | 45.6628 |
| 0.2900 | 99.9776 | 99.9776 | 49.8262 |
| 0.3000 | 99.9816 | 99.9817 | 61.3741 |
| 0.3100 | 99.9849 | 99.9852 | 73.4072 |
| 0.3200 | 99.9875 | 99.9878 | 87.6728 |
| 0.3300 | 99.9897 | 99.9899 | 112.9891 |
| 0.3400 | 99.9914 | 99.9917 | 121.8812 |
| 0.3500 | 99.9929 | 99.9932 | 134.6083 |

The second term satisfies

$$\sum_{i=1}^{n} \xi_i \mathbf{1}_{\{\nu_i=0\}} \le \sum_{i=1}^{\kappa_\infty} \xi_i \mathbf{1}_{\{\nu_i=0\}} < \infty, \tag{6.4}$$

almost surely, because the number of reinitialization remains finite and converges to zero if normalized by $\phi(W_n)$. We can therefore just drop the second term in (6.3) and assume that

$$\theta_n \in \mathcal{K}_{\kappa_n}$$

regardless whether $\nu_n \neq 0$ holds or not. Let $\mathcal{K}$ be an arbitrary compact set. Assumption (4.33) implies

$$w_x(\theta_{n-1})^{q-1}|f(x)|^q \mathbf{1}_{\{\theta\in\mathcal{K}\}} \le C(\mathrm{diam}(\mathcal{K}))^{q-1} W(x)^{q-1}|f(x)|^q. \tag{6.5}$$

By Hölder's inequality, we have for $q < p$ and all $\theta \in \mathcal{K}$,

$$\begin{aligned} m_{q,f}(\theta)\mathbf{1}_{\{\theta\in\mathcal{K}\}} &= \mathbb{E}_{\theta_0}[w_X(\theta)^{q-1}|f(X)|^q\mathbf{1}_{\{\theta\in\mathcal{K}\}}] \\ &\le \|f\|^q_{\theta_0,p} C(\operatorname{diam}(\mathcal{K}))^{q-1}\mathbb{E}_{\theta_0}[W(X)^{p(q-1)/(p-q)}]^{(p-q)/q} \qquad (6.6) \\ &< \infty \end{aligned}$$

as long as $q$ satisfies the condition

$$q < \frac{p(p^*+1)}{p+p^*}. \qquad (6.7)$$

Note that (6.7) implies also $q < p$ because $\frac{p(p^*+1)}{p+p^*} \le p$ for $p \ge 1$. Condition (4.63) implies

$$m_{q,f}(\theta) < \infty \qquad \forall 1 \le q \le 4. \qquad (6.8)$$

Let $\mathcal{K}$ be a compact neighborhood of $\theta^*$. Lebesgue's theorem together with the continuity condition (4.64) and the upper bound (6.5), which is integrable by (6.6), shows that

$$m_{f,q} : \theta \mapsto m_{f,q}(\theta) \qquad (6.9)$$

is continuous for $q \le 4$. Without loss of generality, we assume from now on that $\mathbb{E}[f(X)] = 0$. By assumptions (4.33) and (4.34), we have for $q < p$ and $a > 1$

$$\begin{aligned} &\mathbb{E}[w_{X_n}(\theta_{n-1})^q|f(X_n)|^q\mathbf{1}_{\{\kappa_\infty=j\}} \mid \mathcal{F}_{i-1}] \\ &\quad= \mathbb{E}_{\theta_0}[w_{X_n}(\theta_{n-1})^{q-1}|f(X_n)|^q\mathbf{1}_{\{\kappa_\infty=j\}}] \\ &\quad\le \|f\|^q_{\theta_0,p}\mathbb{E}_{\theta_0}[w_{X_n}(\theta_{n-1})^{p(q-1)/(p-q)}\mathbf{1}_{\{\kappa_\infty=j\}}]^{(p-q)/p} \\ &\quad\le \|f\|^q_{\theta_0,p}\mathbb{P}(\kappa_\infty = j)^{(p-q)/p1/a}, \\ &\mathbb{E}_{\theta_0}[w_{X_n}(\theta_{n-1})^{p(q-1)/(p-q)a/(a-1)}\mathbf{1}_{\{\kappa_\infty=j\}}]^{(p-q)/p(a-1)/a} \\ &\quad\le \|f\|^q_{\theta_0,p}C(\rho)^{(p-q)/p1/a}e^{k_{p^*}(p-q)/p1/a} \\ &\qquad\times(\rho^{(p-q)/p1/a+p(q-1)/(p-q)m_{p^*}\log_\rho(e)})^j\|W\|^{q-1}_{\theta_0,q_1} < \infty, \end{aligned}$$

if $q_1 \le p^*$ holds for

$$q_1 = \frac{p(q-1)}{p-q}\frac{a}{a-1}. \qquad (6.10)$$

We may choose $a$ arbitrarily large at the expense of increasing the constant in the above estimate. Therefore, $q_1 \le p^*$ holds if

$$\frac{p(q-1)}{p-q} < p^* \qquad (6.11)$$

which is equivalent to (6.7). Next, we choose

$$\rho < e^{-am_{p^*}p(q-1)/(p-q)}, \tag{6.12}$$

such that we can sum over $j = 1, \ldots, \infty$ to obtain

$$\begin{aligned} m_{f,q}(\theta_{n-1}) &= \mathbb{E}[w_{X_n}(\theta_{n-1})^q |f(X_n)|^q \mid \mathcal{F}_{n-1}] \\ &= \sum_{j=1}^{\infty} \mathbb{E}[w_{X_n}(\theta_{n-1})^q |f(X_n)|^q \mathbf{1}_{\{\kappa_\infty = j\}} \mid \mathcal{F}_{n-1}] \\ &< C(\rho, a, p, p^*, q, \|f\|_{\theta_0,p}, \|W\|_{\theta_0,p^*}) \end{aligned}$$

with an upper bound independent of $n$. Assumption (4.63) implies that

$$\sup_n \mathbb{E}[m_{f,q}(\theta_{n-1})] < \infty \qquad \forall 1 \le q \le 4. \tag{6.13}$$

We have

$$\langle M \rangle_n = \sum_{i=1}^{n} (m_{f,2}(\theta_{i-1}) - \mathbb{E}[f(X)]^2). \tag{6.14}$$

Because $\theta \mapsto m_{f,2}(\theta)$ is continuous at $\theta^*$ and $\theta_{i-1} \to \theta^*$ almost surely, we obtain from Cesaro's lemma that

$$\frac{\langle M \rangle_n}{n} = \frac{1}{n} \sum_{i=1}^{n} (m_{f,2}(\theta_{i-1}) - \mathbb{E}[f(X)]^2) \quad \to \quad m_{f,2}(\theta^*) - \mathbb{E}[f(X)]^2 = \sigma_f^2(\theta^*),$$

almost surely. By (6.13) and Lebesgue's dominated convergence theorem,

$$\frac{s_n^2}{n} = \frac{\mathbb{E}[\langle M \rangle_n]}{n} \to \sigma_f^2(\theta^*). \tag{6.15}$$

Set

$$\bar{M}_n = \sum_{i=1}^{n} (\xi_i^2 - \mathbb{E}[\xi_i^2 \mid \mathcal{F}_{i-1}]).$$

By (6.13) $\bar{M}_n$ is a square integrable martingale because

$$\begin{aligned} \mathbb{E}[\Delta \bar{M}_i^2 \mid \mathcal{F}_{i-1}] &= \mathbb{E}[(\xi_i^2 - \mathbb{E}[\xi_i^2 \mid \mathcal{F}_{i-1}])^2 \mid \mathcal{F}_{i-1}] \\ &\le 8(m_{f,4}(\theta_{i-1}) + \mathbb{E}[f(X)]^4). \end{aligned}$$

More precisely,

$$\begin{aligned} \mathbb{E}[\Delta \bar{M}_i^2 \mid \mathcal{F}_{i-1}] = m_{f,4}(\theta_{i-1}) &- m_{f,2}(\theta_{i-1})^2 - 4m_{f,3}(\theta_{i-1})\mathbb{E}[f(X)] \\ &+ 8m_{f,2}(\theta_{i-1})\mathbb{E}[f(X)]^2 - 4\mathbb{E}[f(X)]^4. \end{aligned}$$

The continuity of $\theta \mapsto m_{f,p}(\theta)$ in $\theta^*$ for $1 \le p \le 4$ and Cesaro's lemma imply

$$\begin{aligned} &\frac{1}{n}\sum_{i=1}^{n} \mathbb{E}[\Delta \bar{M}_i^2 \mid \mathcal{F}_{i-1}] \\ &\quad \to m_{f,4}(\theta^*) - m_{f,2}(\theta^*)^2 - 4m_{f,3}(\theta^*)\mathbb{E}[f(X)] \\ &\qquad + 8m_{f,2}(\theta^*)\mathbb{E}[f(X)]^2 - 4\mathbb{E}[f(X)]^4, \end{aligned} \tag{6.16}$$

almost surely, which together with (6.15) implies

$$\lim_{n\to\infty} s_n^{-2}\sum_{i=1}^{n}(\xi_i^2 - \mathbb{E}[\xi_i^2 \mid \mathcal{F}_{i-1}]) = \lim_{n\to\infty}\frac{1}{n}\sum_{i=1}^{n}(\xi_i^2 - \mathbb{E}[\xi_i^2 \mid \mathcal{F}_{i-1}]) = 0, \tag{6.17}$$

almost surely. Therefore,

$$\lim_{n\to\infty}\frac{[M]_n}{\langle M\rangle_n} = 1 + \lim_{n\to\infty}\left(\frac{\langle M\rangle_n}{n}\right)^{-1}\frac{1}{n}\sum_{i=1}^{n}(\xi_i^2 - \mathbb{E}[\xi_i^2 \mid \mathcal{F}_{i-1}]) = 1, \tag{6.18}$$

almost surely. To apply Corollary 4.2 in [15], we need to verify the three conditions:

$$s_n^{-2}[M]_n \to \eta^2 > 0 \qquad \text{almost surely;} \tag{6.19}$$

$$\forall \varepsilon > 0 \qquad \sum_{n=1}^{\infty} s_n^{-1}\mathbb{E}[|\xi_n|\mathbf{1}_{\{|\xi_n|>\varepsilon s_n\}}] < \infty; \tag{6.20}$$

$$\exists \delta > 0 \qquad \sum_{n=1}^{\infty} s_n^{-4}\mathbb{E}[|\xi_n|^4\mathbf{1}_{\{|\xi_n|\le\delta s_n\}}] < \infty. \tag{6.21}$$

Condition (6.19) holds because

$$\begin{aligned} \lim_{n\to\infty}\frac{[M]_n}{s_n^2} &= \lim_{n\to\infty}\left(\frac{s_n^2}{n}\right)^{-1}\frac{\langle M\rangle_n}{n} \\ &\quad + \lim_{n\to\infty} s_n^{-2}\sum_{i=1}^{n}(\xi_i^2 - \mathbb{E}[\xi_i^2 \mid \mathcal{F}_{i-1}]) \\ &= 1, \end{aligned} \tag{6.22}$$

almost surely, as a consequence of (6.15) and (6.17). By (6.15), we may replace $s_n^2$ by $n$ for the verification of conditions (6.20) and (6.21).

Let $1 < a < 2$. We first approach (6.20). From Hölder's and Chebyshev's inequalities, we have

$$\begin{aligned} &\sqrt{n}^{-1}\mathbb{E}[|\xi_n|\mathbf{1}_{\{|\xi_n|>\varepsilon\sqrt{n}\}}] \\ &\quad \le \sqrt{n}^{-1}\mathbb{E}[|\xi_n|^{2a}]^{1/(2a)}\mathbb{P}(|\xi_n| > \varepsilon\sqrt{n})^{1-1/(2a)} \end{aligned}$$

(6.23)
$$\leq \sqrt{n}^{-1}\mathbb{E}[|\xi_n|^{2a}]^{1/(2a)}\mathbb{E}[|\xi_n|^{2a}]^{1-1/(2a)}\left(\frac{1}{\varepsilon\sqrt{n}}\right)^{2a-1}$$
$$\leq \varepsilon^{1-2a}\mathbb{E}[|\xi_n|^{2a}]n^{-a}$$

for every fixed $\varepsilon > 0$. Therefore,

$$\sum_{n=1}^{\infty}\sqrt{n}^{-1}\mathbb{E}[|\xi_n|\mathbf{1}_{\{|\xi_n|>\varepsilon\sqrt{n}\}}] \leq \varepsilon^{1-2a}\sum_{n=1}^{\infty}\mathbb{E}[|\xi_n|^{2a}]n^{-a} < \infty. \tag{6.24}$$

This last equation implies condition (6.20). To check condition (6.21), note that

$$\sum_{n=1}^{\infty} n^{-2}\mathbb{E}[|\xi_n|^4\mathbf{1}_{\{|\xi_n|\leq\delta\sqrt{n}\}}] \leq \sum_{n=1}^{\infty} n^{-2}\mathbb{E}[(\delta\sqrt{n})^{4-2a}|\xi_n|^{2a}\mathbf{1}_{\{|\xi_n|\leq\delta\sqrt{n}\}}]$$
(6.25)
$$\leq \delta^{4-2a}\sum_{n=1}^{\infty} n^{-a}\mathbb{E}[|\xi_n|^{2a}] < \infty.$$

The sums (6.24) and (6.25) are finite because

$$\mathbb{E}[|\xi_n|^{2a} \mid \mathcal{F}_{n-1}] \leq 2^{2a-1}(m_{f,2a}(\theta_n) + E[f(X)]^{2a}), \tag{6.26}$$

and $\sup_n \mathbb{E}[m_{f,2a}(\theta_n)] < \infty$, as shown in (6.13).

6.2. *Proof of Theorem 4.2.* Under the assumptions of Theorem 4.2, the boundedness of the functions

$$f_y = \mathbf{1}_{(-\infty,y]} \circ \Psi, \qquad 1 - f_y, \qquad y \in \mathbb{R} \tag{6.27}$$

allows us to apply the law of iterated logarithm (Theorem 4.4). We verify the convergence statement by proving that

$$\mathbb{P}(q_{n,w,\nu}(\alpha) \leq q_\alpha - \delta \text{ i.o.}) = 0 \qquad \forall\delta > 0, \tag{6.28}$$

and

$$\mathbb{P}(q_{n,w,\nu}(\alpha) > q_\alpha \text{ i.o.}) = 0, \tag{6.29}$$

where i.o. stands for infinitely often and is defined as

$$A_n \text{ i.o.} = \limsup_n A_n = \bigcap_{n=1}^{\infty}\bigcup_{k=n}^{\infty} A_k. \tag{6.30}$$

Let $F(y) = \mathbb{P}(Y \leq y)$ denote the distribution function of $Y = \Psi(X)$. We first analyze the estimator $q_{n,w,\nu}(\alpha)$. Define

$$A_n(\delta) = \{q_{n,w,\nu}(\alpha) \leq q_\alpha - \delta\}. \tag{6.31}$$

It follows from (4.26) and Lemma 6.1 that

$$\begin{aligned} A_n(\delta) &= \left\{ \frac{1}{\nu(n)\sum_i w_i} \sum_i w_i \mathbf{1}_{\{Y_i \le q_\alpha - \delta\}} \ge \alpha \right\} \\ &= \left\{ \sum_i (w_i \mathbf{1}_{\{Y_i \le q_\alpha - \delta\}} - F(q_\alpha - \delta)) \ge \nu(n)\alpha \sum_i w_i - nF(q_\alpha - \delta) \right\}. \end{aligned} \tag{6.32}$$

Let

$$W_n(\eta) = \left\{ \left| \sum_i (w_i - 1) \right| \le (1+\eta)\phi(nv) \right\}. \tag{6.33}$$

We consider

$$A_n(\delta) \subset A_n(\delta) \cap W_n(\eta) \cup \complement W_n(\eta). \tag{6.34}$$

Then

$$\begin{aligned} &A_n(\delta) \cap W_n(\eta) \\ &\quad \subset \left\{ \sum_i (w_i \mathbf{1}_{\{Y_i \le q_\alpha - \delta\}} - F(q_\alpha - \delta)) \right. \\ &\qquad \left. \ge \nu(n)\alpha(n - (1+\eta)\phi(nv)) - nF(q_\alpha - \delta) \right\}. \end{aligned} \tag{6.35}$$

Similarly, we have

$$\begin{aligned} B_n &= \{ q_{n,w,\nu}(\alpha) > q_\alpha \} \\ &= \left\{ \sum_i (w_i \mathbf{1}_{\{Y_i \le q_\alpha\}} - F(q_\alpha)) < \nu(n)\alpha \sum_i w_i - nF(q_\alpha) \right\} \end{aligned} \tag{6.36}$$

and

$$\begin{aligned} &B_n \cap W_n(\eta) \\ &\quad \subset \left\{ \sum_i (w_i \mathbf{1}_{\{Y_i \le q_\alpha\}} - F(q_\alpha)) \right. \\ &\qquad \left. < \nu(n)\alpha(n + (1+\eta)\phi(nv)) - nF(q_\alpha) \right\}. \end{aligned} \tag{6.37}$$

For arbitrary $\eta > 0$, let

$$A_n^{\mathrm{LIL}}(\delta, \eta) = \left\{ \sum_i (w_i \mathbf{1}_{\{Y_i \le q_\alpha - \delta\}} - F(q_\alpha - \delta)) \ge (1+\eta)\phi(nv_{q_\alpha - \delta}) \right\} \tag{6.38}$$

and

$$B_n^{\mathrm{LIL}}(\eta) = \left\{ \sum_i (w_i \mathbf{1}_{\{Y_i \le q_\alpha\}} - F(q_\alpha)) \le -(1+\eta)\phi(nv_\alpha) \right\}. \tag{6.39}$$

Then

$$\begin{aligned} &\frac{1+\eta}{\alpha}\phi(nv_{q_\alpha-\delta}) + \frac{F(q_\alpha-\delta)}{\alpha}n \le \nu(n)(n-(1+\eta)\phi(nv)) \\ &\implies \quad A_n(\delta)\cap W_n(\eta) \subset A_n^{\mathrm{LIL}}(\delta,\eta) \end{aligned} \tag{6.40}$$

and

$$\begin{aligned} &\nu(n)(n+(1+\eta)\phi(nv)) \le \frac{F(q_\alpha)}{\alpha}n - \frac{1+\eta}{\alpha}\phi(nv_\alpha) \\ &\implies \quad B_n\cap W_n(\eta) \subset B_n^{\mathrm{LIL}}(\delta,\eta). \end{aligned} \tag{6.41}$$

Recall that

$$\limsup_n (A_n \cup B_n) = \limsup_n A_n \cup \limsup_n B_n. \tag{6.42}$$

Hence, if (6.40) is satisfied, we have

$$\begin{aligned} \mathbb{P}(A_n(\delta)\text{ i.o.}) &\le \mathbb{P}(A_n(\delta)\cap W_n(\eta)\text{ i.o.}) + \mathbb{P}(\complement W_n(\eta)\text{ i.o.}) \\ &\le \mathbb{P}(A_n^{\mathrm{LIL}}(\delta,\eta)\text{ i.o.}) + \mathbb{P}(\complement W_n(\eta)\text{ i.o.}). \end{aligned} \tag{6.43}$$

From the law of iterated logarithm in Theorem 4.4, we know that

$$\mathbb{P}(A_n^{\mathrm{LIL}}(\delta,\eta)\text{ i.o.}) = 0, \qquad \mathbb{P}(\complement W_n(\eta)\text{ i.o.}) = 0. \tag{6.44}$$

Therefore, $\mathbb{P}(A_n(\delta)\text{ i.o.}) = 0$ for all $\delta > 0$. In the same way, we obtain

$$\mathbb{P}(B_n\text{ i.o.}) = 0. \tag{6.45}$$

To verify that condition (4.36) implies (6.41) and (6.41), it is sufficient to note that $F(q_\alpha) \ge \alpha$. Because $F(q_\alpha - \delta) < \alpha$, it follows that

$$\frac{1+\eta}{\alpha}\phi(nv_{q_\alpha-\delta}) + \frac{F(q_\alpha-\delta)}{\alpha}n \le n - kn^{1/2+\gamma}$$

for $n$ large enough and for all $\delta > 0$.

The convergence proof for $q^r_{n,w,\nu}(\alpha)$ is slightly simpler. From (4.29) and Lemma 6.1, we get

$$\begin{aligned} A_n^r(\delta) &= \{q^r_{n,w,\nu}(\alpha) \le q_\alpha - \delta\} \\ &= \left\{ 1 - \frac{1}{\nu(n)}\sum_i w_i \mathbf{1}_{\{Y_i > q_\alpha - \delta\}} \ge \alpha \right\} \\ &= \left\{ \sum_i (w_i \mathbf{1}_{\{Y_i > q_\alpha-\delta\}} - (1 - F(q_\alpha-\delta))) \right. \\ &\qquad \left. \le \nu(n)(1-\alpha) - n(1-F(q_\alpha-\delta)) \right\} \end{aligned} \tag{6.46}$$

and

$$
\begin{aligned}
B_n^r &= \{q_{n,w,\nu}^r(\alpha) > q_\alpha\} \\
&= \left\{ \sum_i (w_i \mathbf{1}_{\{Y_i > q_\alpha\}} - (1 - F(q_\alpha))) \right. \\
&\qquad \left. > \nu(n)(1-\alpha) - n(1 - F(q_\alpha)) \right\}.
\end{aligned} \tag{6.47}
$$

For arbitrary $\eta > 0$, let

$$
\begin{aligned}
A_n^{r,\mathrm{LIL}}(\delta,\eta) = \left\{ \sum_i (w_i \mathbf{1}_{\{Y_i > q_\alpha - \delta\}} - (1 - F(q_\alpha - \delta))) \right. \\
\left. \le -(1+\eta)\phi(n v_{q_\alpha - \delta}) \right\}
\end{aligned} \tag{6.48}
$$

and

$$
B_n^{r,\mathrm{LIL}}(\eta) = \left\{ \sum_i (w_i \mathbf{1}_{\{Y_i > q_\alpha\}} - (1 - F(q_\alpha))) \ge (1+\eta)\phi(n v_\alpha) \right\}. \tag{6.49}
$$

We have

$$
\begin{aligned}
\nu(n) \le \frac{1 - F(q_\alpha - \delta)}{1-\alpha} n - \frac{1+\eta}{1-\alpha}\phi(n v_{q_\alpha - \delta}) \\
\Longrightarrow \quad A_n^r(\delta) \subset A_n^{r,\mathrm{LIL}}(\delta,\eta),
\end{aligned} \tag{6.50}
$$

and

$$
\frac{1+\eta}{1-\alpha}\phi(n v_\alpha) + \frac{1 - F(q_\alpha)}{1-\alpha} n \le \nu(n) \quad \Longrightarrow \quad B_n^r \subset B_n^{r,\mathrm{LIL}}(\eta). \tag{6.51}
$$

By the law of iterated logarithm (Theorem 4.4), we obtain

$$
\mathbb{P}(A_n^{r,\mathrm{LIL}}(\delta,\eta) \text{ i.o.}) = 0, \mathbb{P}(B_n^{r,\mathrm{LIL}}(\eta) \text{ i.o.}) = 0. \tag{6.52}
$$

Therefore, conditions (6.50) and (6.51) are sufficient to guarantee (6.28) and (6.29) for $q_{n,w,\nu}^r$. Because $1 - F(q_\alpha) \le 1 - \alpha$ and $1 - F(q_\alpha - \delta) > 1 - \alpha$ for all $\delta > 0$, condition (4.38) is sufficient for (6.50) and (6.51).

In a completely analogous manner, we obtain

$$
\begin{aligned}
A_n^l(\delta) &= \{q_{n,w,\nu}^l(\alpha) \le q_\alpha - \delta\} \\
&= \left\{ \sum_i (w_i \mathbf{1}_{\{Y_i \le q_\alpha - \delta\}} - F(q_\alpha - \delta)) \ge \nu(n)\alpha - nF(q_\alpha - \delta) \right\}
\end{aligned} \tag{6.53}
$$

and

$$B_n^l = \{q_{n,w,\nu}^l(\alpha) > q_\alpha\} \\ = \Big\{\sum_i (w_i \mathbf{1}_{\{Y_i \le q_\alpha\}} - F(q_\alpha)) < \nu(n)\alpha - nF(q_\alpha)\Big\}. \tag{6.54}$$

This time, let, for $\eta > 0$,

$$A_n^{l,\mathrm{LIL}}(\delta,\eta) = \Big\{\sum_i (w_i \mathbf{1}_{\{Y_i \le q_\alpha - \delta\}} - F(q_\alpha - \delta)) \ge (1+\eta)\phi(nv_{q_\alpha-\delta})\Big\} \tag{6.55}$$

and

$$B_n^{l,\mathrm{LIL}}(\eta) = \Big\{\sum_i (w_i \mathbf{1}_{\{Y_i \le q_\alpha\}} - F(q_\alpha)) \le -(1+\eta)\phi(nv_\alpha)\Big\}. \tag{6.56}$$

We have

$$\frac{1+\eta}{\alpha}\phi(nv_{q_\alpha-\delta}) + \frac{F(q_\alpha-\delta)}{\alpha}n \le \nu(n) \quad \Longrightarrow \quad A_n^l(\delta) \subset A_n^{l,\mathrm{LIL}}(\delta,\eta) \tag{6.57}$$

and

$$\nu(n) \le \frac{F(q_\alpha)}{\alpha}n - \frac{1+\eta}{\alpha}\phi(nv_\alpha) \quad \Longrightarrow \quad B_n^l \subset B_n^{l,\mathrm{LIL}}(\eta). \tag{6.58}$$

By the law of iterated logarithm, equations (6.57) and (6.58) are a sufficient condition to guarantee (6.28) and (6.29) for $q_{n,w,\nu}^l$. Similarly as above, (4.40) is sufficient for (6.57) and (6.58). This proves Theorem 4.2.

6.3. *Proof of Theorem 4.1.* We again apply the law of iterated logarithm 4.4. We only prove the result for $q_{n,w,\nu}(\alpha)$. The other estimators are treated analogously. Because $F$ is increasing in $q_\alpha$, it follows that $F^{\leftarrow}$ is continuous in $\alpha$. It is sufficient to prove for any $\delta > 0$ that

$$\mathbb{P}(q_{n,w,\nu}(\alpha) \le q_\alpha - \delta \text{ i.o.}) = 0 \tag{6.59}$$

and

$$\mathbb{P}(q_{n,w,\nu}(\alpha) > q_\alpha + \delta \text{ i.o.}) = 0. \tag{6.60}$$

For $\nu(n) \equiv 1$, we obtain from (6.40),

$$\frac{1+\eta}{\alpha}\phi(nv_{q_\alpha-\delta}) + \frac{F(q_\alpha-\delta)}{\alpha}n + (1+\eta)\phi(nv) \le n \\ \Longrightarrow \quad A_n(\delta) \cap W_n(\eta) \subset A_n^{\mathrm{LIL}}(\delta,\eta). \tag{6.61}$$

If we define

$$B_n(\delta) = \{q_{n,w,\nu}(\alpha) > q_\alpha + \delta\}$$
$$= \left\{\sum_i (w_i \mathbf{1}_{\{Y_i \le q_\alpha + \delta\}} - F(q_\alpha + \delta)) < \alpha \sum_i w_i - nF(q_\alpha + \delta)\right\}, \tag{6.62}$$

we deduce

$$n \le \frac{F(q_\alpha + \delta)}{\alpha} n - (1+\eta)\phi(nv) - \frac{1+\eta}{\alpha}\phi(nv_{q_\alpha+\delta})$$
$$\implies \quad B_n(\delta) \cap W_n(\eta) \subset B_n^{\mathrm{LIL}}(\delta, \eta). \tag{6.63}$$

For any $\delta > 0$, $F(q_\alpha - \delta) < \alpha$ and $F(q_\alpha + \delta) > \alpha$. Therefore, if $n$ is large enough, conditions (6.61) and (6.63) are satisfied. We conclude as in the proof of Theorem 4.2.

6.4. *Proof of Proposition 4.1.* Let $\mathcal{K}$ be a compact subset of $\mathcal{W}$. We apply Markov's and Burkholder's inequality,

$$\mathbb{P}\Big(\max_{k \le n} \|S_{1,k}(\boldsymbol{\gamma}, \boldsymbol{\epsilon}, \mathcal{K})\| > \delta\Big)$$
$$\le \frac{B_p}{\delta^p}\mathbb{E}\left[\left(\mathbf{1}_{\{n \le \sigma(\mathcal{K},\boldsymbol{\epsilon})\}} \sum_{k=1}^n \gamma_k^2 \|H(X_k, \theta_{k-1}) - h(\theta_{k-1})\|^2\right)^{p/2}\right]$$
$$\le \frac{2^p B_p}{\delta^p}\left(\sum_{k=1}^n \gamma_k^2 \mathbb{E}[\mathbf{1}_{\{k-1 \le n \le \sigma(\mathcal{K},\boldsymbol{\epsilon})\}}(\|H(X_k, \theta_{k-1})\|^p + \|h(\theta_{k-1})\|^p)]^{2/p}\right)^{p/2},$$

where $B_p$ is a universal constant only depending on $p$. To estimate

$$\mathbb{E}[\mathbf{1}_{\{k-1 \le n \le \sigma(\mathcal{K},\boldsymbol{\epsilon})\}}(\|H(X_k, \theta_{k-1})\|^p)]^{2/p},$$

note that by our assumptions

$$\mathbb{E}[\mathbf{1}_{\{k-1 \le n \le \sigma(\mathcal{K},\boldsymbol{\epsilon})\}}\|H(X_k, \theta_{k-1})\|^p]^{2/p}$$
$$= \mathbb{E}[\mathbf{1}_{\{k-1 \le \sigma(\mathcal{K},\boldsymbol{\epsilon})\}}\mathbb{E}_{\theta_{k-1}}[\|H(X_k, \theta_{k-1})\|^p]]^{2/p}$$
$$= \mathbb{E}\left[\mathbf{1}_{\{k-1 \le \sigma(\mathcal{K},\boldsymbol{\epsilon})\}}\mathbb{E}_{\theta_0}\left[\frac{\|H(X_k, \theta_{k-1})\|^p}{w_{\theta_{k-1}}(X_k) W^p(X_k)} W^p(X_k)\right]\right]^{2/p}$$
$$= C_{\mathcal{K}}^2 \mathbb{E}[\mathbf{1}_{\{k-1 \le \sigma(\mathcal{K},\boldsymbol{\epsilon})\}}\mathbb{E}_{\theta_0}[W^p(X_k)]]^{2/p} \le C_{\mathcal{K}}^2 \|W\|_{\theta_0,p}^2,$$

where the constant $C_{\mathcal{K}}$ comes from assumption (4.54). Because $h$ is continuous and $\mathcal{K}$ compact, $\mathbf{1}_{\{k-1\le\sigma(\mathcal{K},\boldsymbol{\epsilon})\}}\|h(\theta_{k-1})\|$ is bounded as well. Therefore, we arrive at the estimate

$$\mathbb{P}\Big(\max_{k\le n}\|S_{1,k}(\boldsymbol{\gamma},\boldsymbol{\epsilon},\mathcal{K})\|>\delta\Big)\le C\frac{1}{\delta^p}\left(\sum_{k=1}^{n}\gamma_k^2\right)^{p/2}.$$

The bound

$$\mathbb{P}^{\boldsymbol{\gamma}}_{\Phi(x,\theta)}(\nu(\boldsymbol{\epsilon})<\mathcal{K})\le C\sum_{k=1}^{n}\left(\frac{\gamma_k}{\epsilon_k}\right)^p$$

is derived similarly as in the proof of Proposition 5.2 in [1] .

**Acknowledgments.** We would like to thank Michael Wolf and seminar participants at the IBM Research Lab, Rüschlikon, and the Deutsche Bundesbank in Frankfurt for their comments. Part of this research has been carried out within the National Center of Competence in Research "Financial Valuation and Risk Management" (NCCR FINRISK). The NCCR FINRISK is a research program supported by the Swiss National Science Foundation.

QuantCatalyst
8000 Zurich
Switzerland
E-mail: daniel.egloff@hispeed.ch

Swiss Banking Institute
University of Zurich
Plattenstrasse 14
8032 Zurich
Switzerland
E-mail: leippold@isb.uzh.ch